\documentclass[preprint,12pt]{elsarticle}
\usepackage{amsmath,amssymb,amsfonts}
\usepackage{algorithmic}
\usepackage{graphicx}
\usepackage{textcomp}

\usepackage[scr = dutchcal]{mathalpha}
\usepackage{graphicx}
\usepackage{textcomp}
\usepackage{stmaryrd}
\usepackage{tabularx}
\usepackage{multirow}
\usepackage{float}
\usepackage{color}

\usepackage[font=footnotesize]{caption}
\usepackage{subcaption}

\newcommand{\mbf}[1]{\mathbf{ #1}}
\newcommand{\tbf}[1]{\textbf{#1}}
\newcommand{\mbs}[1]{\boldsymbol{#1}}

\newcommand{\mcl}[1]{\mathcal{#1}}
\newcommand{\mscr}[1]{\mathscr{#1}}

\newcommand{\tnf}[1]{\textnormal{#1}}
\newcommand{\R}{\mathbb{R}}
\newcommand{\N}{\mathbb{N}}
\newcommand{\Z}{\mathbb{Z}}
\newcommand{\norm}[1]{\left\lVert{#1}\right\rVert}

\newcommand{\bmat}[1]{\begin{bmatrix}#1\end{bmatrix}}

\newtheorem{thm}{Theorem}

\newtheorem{lem}[thm]{Lemma}
\newtheorem{prop}[thm]{Proposition}
\newtheorem{cor}[thm]{Corollary}
\newtheorem{remark}[thm]{Remark}
\newtheorem{example}[thm]{Example}
\newproof{pf}{Proof}
\newproof{pot}{Proof of Theorem~\ref{thm:fundamental_expansion_ND}}
\newtheorem{thm1}{Theorem}

\let\bl\bigl
\let\bbl\Bigl
\let\bbbl\biggl
\let\bbbbl\Biggl
\let\br\bigr
\let\bbr\Bigr
\let\bbbr\biggr
\let\bbbbr\Biggr

\begin{document}

\begin{frontmatter}
	
\title{\LARGE
		Constructive Representation of Functions in $N$-Dimensional Sobolev Space%
\tnoteref{t1,t2}}
\tnotetext[t1]{This work was supported by National Science Foundation grant CMMI-1935453.}
\author[1]{Declan S. Jagt\corref{cor1}%
	}

\author[1]{Matthew M. Peet%
	}



\affiliation[1]{organization={SEMTE Arizona State University},
	city={Tempe},
	postcode={85287}, 
	state={Arizona},
	country={USA}}

\cortext[cor1]{Corresponding author}	
	
\begin{abstract}
	A new representation is proposed for functions in a Sobolev space with dominating mixed smoothness on an $N$-dimensional hyperrectangle. In particular, it is shown that these functions can be expressed in terms of their highest-order mixed derivative, as well as their lower-order derivatives evaluated along suitable boundaries of the domain.	The proposed expansion is proven to be invertible, uniquely identifying any function in the Sobolev space with its derivatives and boundary values. Since these boundary values are either finite-dimensional, or exist in the space of square-integrable functions, this offers a bijective relation between the Sobolev space and $L_{2}$. Using this bijection, it is shown how approximation of functions in Sobolev space can be performed in the less restrictive space $L_{2}$, reconstructing such an approximation of the function from an $L_{2}$-optimal projection of its boundary values and highest-order derivative. This approximation method is presented using a basis of Legendre polynomials and a basis of step functions,	and results using both bases are demonstrated to exhibit better convergence behavior than a direct projection approach for two numerical examples.
\end{abstract}

%
%
%
%
%

\end{frontmatter}



\section{Introduction}\label{sec:introduction}

We consider the basic question of the relation between a function and its highest-order partial derivative. In particular, given a differentiable function $\mbf{u}$, can we uniquely represent $\mbf{u}$ by its highest-order derivative and a set of independent boundary values? Conversely, given these boundary values, can we reconstruct the associated function $\mbf{u}$ from its highest-order derivative?

To illustrate, for a first-order differentiable function in a single variable, an answer to both of these questions is given by the fundamental theorem of calculus, or Newton-Leibniz formula, which states that a differentiable function $\mbf{u}:[a,b]\to \R$ can be represented by its partial derivative $\partial_{s} \mbf{u}$ and the value of $\mbf{u}$ at a single point in its domain, i.e. $\mbf{u}(s)=\mbf{u}(a)+\int_{a}^{s}\partial_{s}\mbf{u}(\theta)d\theta$. Moreover, given any $c\in\R$ and $\mbf{v}\in L_2[a,b]$ such that $\mbf{u}(s)=c+\int_{a}^{s}\mbf{v}(\theta)d\theta$, it follows that $\mbf{u}(a)=c$ and $\partial_{s}\mbf{u}=\mbf{v}$. If the function is $n$th-order differentiable, this result may be generalized using Cauchy's formula for repeated integration, uniquely expressing $\mbf{u}$ in terms of its highest-order derivative $\partial_{s}^{n}\mbf{u}$ and suitable boundary values (see e.g. Theorem~7.6 in~\cite{apostol1967calculus}).

Now, for functions in $N\in\N$ variables,
consider the Sobolev space $S_{2}^{\delta}$ of functions with bounded mixed derivatives up to order $\delta\in\N^{N}$, so that $D^{\alpha}\mbf{u}=\frac{\partial^{\alpha_{1}}}{\partial s_{1}^{\alpha_{1}}}\cdots\frac{\partial^{\alpha_{N}}}{\partial s_{N}^{\alpha_{N}}}\mbf{u}\in L_2$ for every $\alpha\in\N^{N}$ with $\alpha_{i}\leq \delta_{i}$ for all $i$.
Using the theory of Green's functions, it has been proven that for a variety of differential operators $D$, any $\mbf{u}\in S_{2}^{\delta}$ can be expressed in terms of $D\mbf{u}$ using a suitable integral operator, recovering $\mbf{u}$ from e.g. its gradient $\nabla\mbf{u}$ or Laplacian $\nabla^2\mbf{u}$ if $\delta_{i}\geq 1$ or $\delta_{i}\geq 2$ for all $i$, respectively~\cite{polyanin2001pdes}.
In addition, a generalization of the Newton-Leibniz formula has recently been proposed for the case that $\delta=(1,\hdots,1)\in\N^{N}$, proving that any $\mbf{u}\in S_{2}^{(1,\hdots,1)}$ can be expressed directly in terms of its mixed derivative $D^{(1,\hdots,1)}\mbf{u}\in L_{2}$~\cite{abdulla2023SobolevEmbedding}.
Unfortunately, however, both the results on Green's functions and the generalized Newton-Leibniz formula are still limited to very particular differential operators (e.g. $\nabla$ and $D^{(1,\hdots,1)}$), and an explicit relation between functions $\mbf{u}\in S_{2}^{\delta}$ and their highest-order derivatives $D^{\delta}\mbf{u}$ has not been presented for general $\delta\in\N^{N}$.

In this paper, an explicit relation is proposed between functions $\mbf{u}\in S_{2}^{\delta}[\Omega]$ on a hyperrectangle $\Omega:=\prod_{i=1}^{N}[a_{i},b_{i}]$, and their 
mixed derivatives $D^{\alpha}\mbf{u}$ for $\alpha_{i}\leq\delta_{i}$ for all $i$.
In particular, the main result is as follows.

\begin{thm}\label{thm:fundamental_expansion_ND}
	Suppose $\mbf{u}\in S_{2}^{\delta}[\Omega]$ for $\delta\in\N^{N}$ and $\Omega:=\prod_{i=1}^{N}[a_{i},b_{i}]\subseteq\R^{N}$. Then
	\begin{equation*}
		\mbf{u}(s)=\sum_{0\leq \alpha\leq\delta}\bl(\mcl{G}^{\delta}_{\alpha}B^{\alpha-\delta}D^{\alpha}\mbf{u}\br)(s),\qquad s\in\Omega,
	\end{equation*}
	where $B^{\beta}:=\prod_{i=1}^{N}\tnf{b}_{i}^{\beta_{i}}$ and $\mcl{G}^{\delta}_{\alpha}:=\prod_{i=1}^{N}\mscr{g}^{\delta_{i}}_{i,\alpha_{i}}$, with
	\begin{equation*}
		(\tnf{b}_{i}^{\beta_{i}}\mbf{u})(s)=\begin{cases}
			\mbf{u}(s_{1},\hdots,a_{i},\hdots,s_{N}),	&	\beta_{i}<0,	\\
			\mbf{u}(s_{1},\hdots,s_{i},\hdots,s_{N}),	&	\beta_{i}=0,	
		\end{cases}
	\end{equation*}
	\begin{equation*}
		(\mscr{g}^{\delta_{i}}_{i,\alpha_{i}} \mbf{u})(s)=\begin{cases}
			\mbf{p}_{\alpha_{i}}(s_{i}-a_{i})\mbf{u}(s),	& \alpha_{i}<\delta_{i},\\
			\int_{a_{i}}^{s_{i}} \mbf{p}_{\alpha_{i}-1}(s_{i}-\theta)\mbf{u}(s)|_{s_{i}=\theta}d \theta & \alpha_{i}=\delta_{i},
		\end{cases}
	\end{equation*}
	where $\mbf{p}_{k}(z):=\frac{z^{k}}{k!}$.
	
	Moreover, for any $\{\mbf{v}^{\alpha}\subseteq L_2[\Gamma^{\alpha}]\mid 0\leq \alpha\leq\delta\}$ on suitable $\Gamma^{\alpha}\subseteq\Omega$, if
	\begin{equation*}
		\mbf{u}(s)=\sum_{0\leq \alpha\leq\delta}\bl(\mcl{G}_{\alpha}^{\delta}\mbf{v}^{\alpha}\br)(s),\qquad s\in\Omega,
	\end{equation*}
	then, $\mbf{v}^{\alpha}=B^{\alpha-\delta}D^{\alpha}\mbf{u}$ for all $0\leq \alpha\leq \delta$.
\end{thm}

Theorem~\ref{thm:fundamental_expansion_ND} shows that any function $\mbf{u}\in S_{2}^{\delta}[\Omega]$ can be represented uniquely by its partial derivatives $D^{\alpha}\mbf{u}$ for $0\leq \alpha\leq\delta$, evaluated at suitable (lower) boundaries of the domain. For example, consider a second-order differentiable function $\mbf{u}\in S_{2}^{(2,1)}[[a,b]\times[c,d]]$ in two variables. Then, the associated boundary values in Theorem~\ref{thm:fundamental_expansion_ND} are given by
\begin{align*}
	B^{-(2,1)}D^{(0,0)}\mbf{u}&=\mbf{u}(a,c),	&	(B^{-(2,0)}D^{(0,1)}\mbf{u})(y)&=\partial_{y}\mbf{u}(a,y),	\\
	B^{-(1,1)}D^{(1,0)}\mbf{u}&=\partial_{x}\mbf{u}(a,c),	&	(B^{-(1,0)}D^{(1,1)}\mbf{u})(y)&=\partial_{x}\partial_{y}\mbf{u}(a,y),	\\
	(B^{-(0,1)}D^{(2,0)}\mbf{u})(x)&=\partial_{x}^2\mbf{u}(x,c),	&	(B^{(0,0)}D^{(2,1)}\mbf{u})(x,y)&=\partial_{x}^2\partial_{y}\mbf{u}(x,y),	
\end{align*}
and Theorem~\ref{thm:fundamental_expansion_ND} implies
\begin{align*}
	&\mbf{u}(x,y) 
	= \mbf{u}(a,c)
	+ (x-a)\partial_{x}\mbf{u}(a,c)   
	+\int_{a}^{x}\!\!(x-\theta)\partial_{x}^2\mbf{u}(\theta,c)d\theta	\nonumber\\ 
	&+\int_{c}^{y}\!\!\partial_{y}\mbf{u}(a,\eta)d\eta 
	+\int_{c}^{y}\!\!(x-a)\partial_{x}\partial_{y}\mbf{u}(a,\eta)d\eta     +\int_{a}^{x}\!\!\int_{c}^{y}(x-\theta)
	\partial_{x}^2\partial_{y}\mbf{u}(\theta,\eta)d\eta d\theta.
\end{align*}

Theorem~\ref{thm:fundamental_expansion_ND} is proven in Section~\ref{sec:main_proof}. In the subsequent section, it is then shown how this result may be applied to perform optimal polynomial approximation of functions $\mbf{u}\in S_{2}^{\delta}[\Omega]$, reconstructing such an approximation from a projection of the derivatives $B^{\alpha-\delta}D^{\alpha}\mbf{u}$ onto a polynomial subspace. This approach has the advantage that, unlike a direct projection of $\mbf{u}$ (using the $L_{2}$ inner product), it accurately captures each of the derivatives $D^{\alpha}\mbf{u}$, thus converging in the Sobolev norm. Moreover, this approach can be adapted to perform projection using a basis of step functions, avoiding undesired oscillatory behavior inherent to smooth basis functions, whilst still capturing each of the derivatives of the function.



\section{Notation}\label{sec:notation}

Denote the set of integers as $\Z$, and that of natural numbers as $\N$, writing $\N_{0}:=\N\cup\{0\}$.
For $N\in\N$ and $\alpha,\beta\in\Z^{N}$, write $\alpha\leq \beta$ if $\alpha_{i}\leq \beta_{i}$ for all $i\in\{1,\hdots,N\}$, and write $\alpha<\beta$ if $\alpha\leq\beta$ and $\alpha\neq\beta$. 
Write $0^{N}$ and $1^{N}$ for the vectors of all zeros and ones in $\N_{0}^{N}$, respectively. 

The focus in this paper will be on functions on a hyperrectangle $\Omega:=\Omega_{1}\times\hdots\Omega_{N}$, where $\Omega_{i}=[a_{i},b_{i}]\subseteq\R$ for $i\in\{1,\hdots,N\}$. 
To distinguish the boundaries and interior of this domain, we introduce the notation
\begin{equation*}
	\Omega^{\beta}:=\Omega_{1}^{\beta_{1}}\times\hdots\times\Omega_{N}^{\beta_{N}},\qquad
	\Omega_{i}^{\beta_{i}}:=\begin{cases}
		\{a_{i}\},	&	\beta_{i}<0,	\\
		(a_{i},b_{i}),	&	\beta_{i}=0,	\\
		\{b_{i}\},	&	\beta_{i}>0,
	\end{cases}
\end{equation*}
for all $\beta\in\Z^{N}$, so that $\Omega=\bigcup_{\beta\in\{-1,0,1\}^{N}}\Omega^{\beta}$. The Hilbert space of (equivalence classes of) square-integrable functions on $\Omega^{\beta}$ is denoted by $L_2[\Omega^{\beta}]$. If $\Omega^{\beta}$ corresponds to a boundary of the domain, i.e. $\Omega_{i}^{\beta_{i}}=\{a_{i}\}$ or $\Omega_{i}^{\beta_{i}}=\{b_{i}\}$ for some $i$, then $L_2[\Omega^{\beta}]$ is identified with the space of square-integrable functions along this boundary, so that e.g. $L_2[[a,b]\times\{c\}]\cong L_2[[a,b]]$ and $L_2[\{b\}\times[c,d]]\cong L_2[[c,d]]$. If $\Omega^{\beta}$ corresponds to a vertex of the hyperrectangle, we let $L_2[\Omega^{\beta}]=\R$, so that e.g. $L_2[\{a\}\times\{d\}]\cong\R$.

For a suitably differentiable function $\mbf{u}\in L_2[\Omega]$, define the partial differential operator $D^{\alpha}$ as
\begin{align*}
	D^{\alpha}\mbf{u}=\partial_{s_{1}}^{\alpha_{1}}\cdots\partial_{s_{N}}^{\alpha_{N}}\mbf{u}, \qquad\text{where}\qquad
	\partial^{\alpha_{i}}_{s_{i}}\mbf{u}=\frac{\partial^{\alpha_{i}}}{\partial s_{i}^{\alpha_{i}}}\mbf{u}.
\end{align*}

\subsection{A Sobolev Space with Bounded Mixed Derivatives up to Order $\delta$}

The focus of this paper is on functions in a Sobolev space with dominating mixed smoothness. In particular, for a given order $\delta\in\N_{0}^{N}$ and a hyperrectangle $\Omega=\prod_{i=1}^{N}[a_{i},b_{i}]$, define
\begin{equation}\label{eq:Sobolev_space}
	S_{2}^{\delta}[\Omega]:=\bl\{\mbf{u}\in L_2[\Omega]\mid D^{\alpha}\mbf{u}\in L_2[\Omega],~\forall \alpha\in\N_{0}^{N}:\alpha\leq\delta\},
\end{equation}
with associated norm $\|\mbf{u}\|_{S_{2}^{\delta}}^2:=\sum_{0^{N}\leq\alpha\leq\delta}\|D^{\alpha}\mbf{u}\|_{L_{2}}^2$, and identifying elements that are equal almost everywhere. Note that this definition is distinct from that of standard Sobolev spaces $W_{2}^{k}[\Omega]$ for $k\in\N$, in that every mixed derivative up to and including $D^{\delta}\mbf{u}$ is required to be bounded. Specifically, where $\mbf{u}\in W_{2}^{k}[\Omega]$ in the standard $k$th-order Sobolev space must satisfy $D^{\alpha}\mbf{u}\in L_2[\Omega]$ only for $\alpha\in\N_{0}^{N}$ for which $\|\alpha\|_{1}=\sum_{i=1}^{N}|\alpha_{i}|\leq k$, a function $\mbf{u}\in S_{2}^{\delta}[\Omega]$ must satisfy $D^{\alpha}\mbf{u}\in L_2[\Omega]$ for all $\alpha\in\N_{0}^{N}$ such that $\alpha\leq\delta$. In this sense, even for $\delta=k\cdot 1^{N}$, the space $S_{2}^{\delta}[\Omega]$ is still a proper subspace of $W_{2}^{k}[\Omega]$, and in fact, is the largest subspace of $W_{2}^{k}[\Omega]$ that is embedded into H\"older space~\cite{abdulla2023SobolevEmbedding}.
Working in this subspace simplifies analysis in the following two crucial manners.
\begin{enumerate}
	\item 
	Functions in $S_{2}^{\delta}[\Omega]$ can be differentiable up to different orders $\delta_{i}\in\N_{0}$ along different directions $s_{i}\in[a_{i},b_{i}]$. In particular, it is clear that if $\mbf{u}\in S_{2}^{\delta}[\Omega]$ for some $\delta\in\N_{0}^{N}$, then $D^{\alpha}\mbf{u}\in S_{2}^{\delta-\alpha}[\Omega]$ for all $0^{N}\leq\alpha\leq\delta$.
	
	\item 
	Functions in $S_{2}^{\delta}[\Omega]$ are sufficienty regular to be evaluated along the boundary of the domain. Specifically, for $\delta\geq 1^{N}$, it has been shown that $S_{2}^{\delta}[\Omega]$ is embedded into a space of H\"older continuous functions~\cite{abdulla2023SobolevEmbedding}. In \ref{appx:SobolevSpace}, it is proven that also for $\delta\leq 1^{N}$, if $\delta_{i}>0$ for some $i\in\{1,\hdots,N\}$, then $\mbf{u}\in S_{2}^{\delta}[\Omega]$ is absolutely continuous along $s_{i}=[a_{i},b_{i}]$, and therefore the boundary values $\mbf{u}|_{s_{i}=a_{i}}$ and $\mbf{u}|_{s_{i}=b_{i}}$ are well-defined.
\end{enumerate}

These two properties of the space $S_{2}^{\delta}[\Omega]$ will be vital for proving Theorem~\ref{thm:fundamental_expansion_ND} in the next section. To compactly state and prove that result, define the boundary operator $B^{\beta}:S_{2}^{\delta}[\Omega]\to L_{2}[\Omega^{\beta}]$ by the restriction $B^{\beta}\mbf{u}=\mbf{u}|_{\Omega^{\beta}}$ to the subdomain $\Omega^{\beta}$, which is well-defined for any $\beta\in\Z^{N}$ such that $|\beta_{i}|\leq|\delta_{i}|$.

\section{Proof of Main Result}\label{sec:main_proof}

In this section, the main theorem is restated and proved.

\begin{thm1}
	Suppose $\mbf{u}\in S_{2}^{\delta}[\Omega]$ for $\delta\in\N_{0}^{N}$ and $\Omega:=\prod_{i=1}^{N}[a_{i},b_{i}]\subseteq\R^{N}$. Then
	\begin{equation}\label{eq:fundamental_expansion_ND}
	\mbf{u}(s)=\sum_{0^{N}\leq \alpha\leq \delta}\bl(\mcl{G}^{\delta}_{\alpha}B^{\alpha-\delta}D^{\alpha}\mbf{u}\br)(s),\qquad s\in\Omega,
	\end{equation}
	where $\mcl{G}^{\delta}_{\alpha}:=\prod_{i=1}^{N}\mscr{g}^{\delta_{i}}_{i,\alpha_{i}}$ with for all $i\in\{1,\hdots,N\}$,
	\begin{equation*}
	(\mscr{g}^{\delta_{i}}_{i,\alpha_{i}} \mbf{u})(s)=\begin{cases}
		\mbf{u}(s),	&	0=\alpha_{i}=\delta_{i},	\\
		\mbf{p}_{\alpha_{i}}(s_{i}\!-\!a_{i})\mbf{u}(s),	& 0\leq \alpha_{i}<\delta_{i},\\
		\int_{a_{i}}^{s_{i}} \mbf{p}_{\alpha_{i}-1}(s_{i}\!-\!\theta)\mbf{u}(s)|_{s_{i}=\theta}d \theta & 0<\alpha_{i}=\delta_{i},
	\end{cases}
	\end{equation*}
	where $\mbf{p}_{k}(z):=\frac{z^{k}}{k!}$.
	Moreover, for any $\{\mbf{v}^{\alpha}\in L_2[\Omega^{\alpha-\delta}]\mid 0^{N}\leq \alpha\leq \delta\}$, if
	\begin{equation}\label{eq:inverse_expansion_ND}
		\mbf{u}(s)=\sum_{0\leq \alpha\leq \delta}\bl(\mcl{G}^{\delta}_{\alpha}\mbf{v}^{\alpha}\br)(s),\qquad s\in\Omega,
	\end{equation}
	then, $\mbf{v}^{\alpha}=B^{\alpha-\delta}D^{\alpha}\mbf{u}$ for all $0^{N}\leq \alpha\leq \delta$.
\end{thm1}

The theorem shows that any function $\mbf{u}\in S_{2}^{\delta}[\Omega]$ can be expressed in terms of each of its partial derivatives $D^{\alpha}\mbf{u}$, evaluated along a suitable lower boundary $\Omega^{\alpha-\delta}$ of the domain. To prove this result, it is first shown how such a function $\mbf{u}$ can be expressed in terms of its derivatives $\partial_{s_{i}}^{\alpha_{i}}\mbf{u}$ with respect to just a single variable $s_{i}$, evaluated at the boundary $s_{i}=a_{i}$. 
To illustrate, consider a function $\mbf{u}\in S_{2}^{(2,1)}[[a,b]\times[c,d]]$ in just two variables. Applying the fundamental theorem of calculus twice, $\mbf{f}$ can be expanded along $x\in[a,b]$ as
\begin{align*}
	\mbf{u}(x,y)&=\mbf{u}(a,y)+\int_{a}^{x}\partial_{x}\mbf{u}(\theta,y)d\theta	\\
	&=\mbf{u}(a,y) +[x-a]\partial_{x}\mbf{u}(a,y) +\int_{a}^{x}\!\!\int_{a}^{\theta}\partial_{x}^2\mbf{u}(\eta,y)d\eta d\theta.
\end{align*}
Using Cauchy's formula for repeated integration, the double integral in this expansion can be expressed as a single integral, yielding a representation of $\mbf{u}$ as
\begin{equation*}
	\mbf{u}(x,y)=\mbf{u}(a,y) +[x-a]\partial_{x}\mbf{u}(a,y) +\int_{a}^{x}[x-\theta]\partial_{x}^2\mbf{u}(\theta,y) d\theta.
\end{equation*}
The following lemma generalizes this result to functions in $N$ variables, differentiable up to arbitrary order $\delta_{i}\in\N_{0}$ with respect to variable $s_{i}\in[a_{i},b_{i}]$, for any $i\in\{1,\hdots,N\}$. In this lemma, $\text{e}_{i}\in\R^{N}$ denotes the $i$th standard Euclidean basis vector in $\R^{N}$ for every $i\in\{1,\hdots,N\}$.

\begin{lem}\label{lem:fundamental_expansion_1D}
	Suppose $\mbf{u}\in S_{2}^{\delta}[\Omega]$ for $\delta\in\N_{0}^{N}$ and $\Omega:=\prod_{i=1}^{N}[a_{i},b_{i}]\subseteq\R^{N}$. Then, for any $i\in\{1,\hdots,N\}$,
	\begin{equation}\label{eq:fundamental_expansion_1D}
		\mbf{u}(s)=\sum_{k=0}^{\delta_{i}}\bl(\mscr{g}^{\delta_{i}}_{i,k}\tnf{b}_{i}^{k-\delta_{i}}\partial_{s_{i}}^{k}\mbf{u}\br)(s),\qquad s\in\Omega,
	\end{equation}
	where $\mscr{g}_{i,k}^{\delta_{i}}$ is as in Theorem~\ref{thm:fundamental_expansion_ND} and $\tnf{b}^{k-\delta_{i}}=B^{(k-\delta_{i})\text{e}_{i}}$, so that
	\begin{align*}
		(\tnf{b}_{i}^{\beta_{i}}\mbf{u})(s)&=\begin{cases}
			\mbf{u}(s_{1},\hdots,a_{i},\hdots,s_{N}),	&	\beta_{i}<0,	\\
			\mbf{u}(s_{1},\hdots,s_{i},\hdots,s_{N}),	&	\beta_{i}=0.	
		\end{cases}
	\end{align*}
	Moreover, for any $\{\mbf{v}^{k}\in L_2[\Omega^{(k-\delta_{i})\text{e}_{i}}]\mid 0\leq k\leq \delta_{i}\}$, if
	\begin{equation}\label{eq:inverse_expansion_1D}
		\mbf{u}(s)=\sum_{k=0}^{\delta_{i}}\bl(\mscr{g}^{\delta_{i}}_{i,k}\mbf{v}^{k}\br)(s),\qquad s\in\Omega,
	\end{equation}
	then, $\mbf{v}^{k}=\tnf{b}_{i}^{k-\delta_{i}}\partial_{s_{i}}^{k}\mbf{u}$ for all $0\leq k\leq \delta_{i}$.
\end{lem}

\begin{pf}
	Without loss of generality, fix $i=1$, re-ordering the variables $s_{i}$ if necessary. Let $\delta\in\N_{0}^{N}$. The result is proven through induction on the order of differentiability $\delta_{i}=\delta_{1}\in\N_{0}$.
	
	\tbf{Base case }$\mbs{\delta_{1}=0:}$
	For the case $\delta_{1}=0$, the result holds trivially, since $\mscr{g}_{1,0}^{0}\tnf{b}_{1}^{0}\partial_{s_{1}}^{0}=I$ is an identity operator.

	\tbf{Induction step }$\mbs{\delta_{1}>0:}$
	To show that the result holds for all $\delta_{1}\in\N_{0}$, fix arbitrary $\tilde{\delta}_{1}\in\N_{0}$. Then, for every $\{\mbf{v}^{k}\in L_2[\Omega^{(k-\tilde{\delta}_{1})\text{e}_{1}}]\mid 0\leq k\leq \tilde{\delta}_{1}\}$,
	\begin{align}\label{eq:fund_int}
		\int_{a_{1}}^{s_{1}}&\bbbbl(
		\sum_{k=0}^{\tilde{\delta}_{1}}\bl(\mscr{g}^{\tilde{\delta}_{1}}_{1,k}\mbf{v}^{k}\br)(\eta) \bbbbr) d\eta	\notag\\
		=&\!\sum_{k=0}^{\tilde{\delta}_{1}-1}\int_{a_{1}}^{s_{1}}\!\!\mbf{p}_{k}(\eta\!-\!a_{1})\mbf{v}^{k}\: d\eta +\!\!\int_{a_{1}}^{s_{1}}\!\!\!\int_{a_{1}}^{\eta}\!\!\mbf{p}_{\tilde{\delta}_{1}-1}(\eta\!-\!\theta)\mbf{v}^{\tilde{\delta}_{1}}(\theta)\: d\theta d\eta	\notag\\
		=&\!\sum_{k=0}^{\tilde{\delta}_{1}-1}\int_{a_{1}}^{s_{1}}\!\!\mbf{p}_{k}(\eta\!-\!a_{1}) d\eta\:\mbf{v}^{k} +\!\!\int_{a_{1}}^{s_{1}}\!\!\!\int_{\theta}^{s_{1}}\!\!\!\mbf{p}_{\tilde{\delta}_{1}-1}(\eta\!-\!\theta)d\eta\:\mbf{v}^{\tilde{\delta}_{1}}(\theta) d\theta 	\notag\\
		=&\!\sum_{k=0}^{\tilde{\delta}_{1}-1}\mbf{p}_{k+1}(s_{1}\!-\!a_{1}) \mbf{v}^{k} +\int_{a_{1}}^{s_{1}}\mbf{p}_{\tilde{\delta}_{1}}(s_{1}\!-\!\theta)\mbf{v}^{\tilde{\delta}_{1}}(\theta) d\theta	
		\hspace*{1.0cm}=\sum_{k=0}^{\tilde{\delta}_{1}}\bl(\mscr{g}^{\tilde{\delta}_{1}+1}_{1,k+1}\mbf{v}^{k}\br)(s_{1}),
	\end{align}
	where the fact that $\int_{y}^{x}\mbf{p}_{k}(\eta-y)d\eta=\mbf{p}_{k+1}(x-y)$ follows by definition of the polynomials $\mbf{p}_{k}$. Now, assume for induction that the lemma holds for $\tilde{\delta}\in\N_{0}^{N}$, and let $\delta=\tilde{\delta}+\text{e}_{1}$ so that $\delta_{1}=\tilde{\delta}_{1}+1$. Fix arbitrary $\mbf{u}\in S_{2}^{\delta}[\Omega]$. To prove that~\eqref{eq:fundamental_expansion_1D} holds for $\mbf{u}$, note that $\partial_{s_{1}}\mbf{u}\in S_{2}^{\delta-\text{e}_{1}}[\Omega]=S_{2}^{\tilde{\delta}}[\Omega]$.
	Applying the fundamental theorem of calculus, invoking the induction hypothesis, and finally using the relation in~\eqref{eq:fund_int}, it follows that
	\begin{align*}
		&\mbf{u}(s_{1})=\mbf{u}(a_{1})+\int_{a_{1}}^{s_{1}}\partial_{s_{1}}\mbf{u}(\eta)d\eta	\\
		&=\mbf{u}(a_{1}) +\int_{a_{1}}^{s_{1}}\bbbbl(\sum_{k=0}^{\tilde{\delta}_{1}}\mscr{g}^{\tilde{\delta}_{1}}_{1,k}\tnf{b}_{1}^{k-\tilde{\delta}_{1}}\partial_{s_{1}}^{k}\bl(\partial_{s_{1}}\mbf{u}\br)\bbbbr)(\eta) d\eta	\\
		&=\mscr{g}_{1,0}^{\delta_{1}}\tnf{b}_{1}^{-\delta_{1}}\partial_{s_{1}}^{0}\mbf{u} +\sum_{k=0}^{\tilde{\delta}_{1}}\bbl(\mscr{g}^{\tilde{\delta}_{1}+1}_{1,k+1}\tnf{b}_{1}^{k-\tilde{\delta}_{1}}\partial_{s_{1}}^{k}\bl(\partial_{s_{1}}\mbf{u}\br)\bbr)(s_{1})	\\
		&=\mscr{g}_{1,0}^{\delta_{1}}\tnf{b}_{1}^{-\delta_{1}}\partial_{s_{1}}^{0}\mbf{u} +\sum_{k=0}^{\delta_{1}-1}\bbl(\mscr{g}^{\delta_{1}}_{1,k+1}\tnf{b}_{1}^{k+1-\delta_{1}}\partial_{s_{1}}^{k+1}\mbf{u}\bbr)(s_{1})	
		=\sum_{k=0}^{\delta_{1}}\bbl(\mscr{g}_{1,k}^{\delta_{1}}\tnf{b}_{1}^{k-\delta_{1}}\partial_{s_{1}}^{k}\mbf{u}\bbr)(s_{1}).
	\end{align*}
	Thus~\eqref{eq:fundamental_expansion_1D} holds. To see that also the implication for~\eqref{eq:inverse_expansion_1D} holds, suppose that for some $\{\mbf{v}^{k}\in L_2[\Omega^{(k-\delta_{i})\text{e}_{i}}]\mid 0\leq k\leq\delta_{1}\}$, we have $\mbf{u}(s_{1})=\sum_{k=0}^{\delta_{1}}\bl(\mscr{g}_{1,k}^{\delta_{1}}\mbf{v}^{k}\br)(s_{1})$. Then, using Eqn.~\eqref{eq:fund_int},
	\begin{equation*}
		\mbf{u}(s_{1})=\mbf{v}^{0} +\sum_{k=0}^{\delta_{1}-1}\bl(\mscr{g}_{1,k+1}^{\delta_{1}}\mbf{v}^{k+1}\br)(s_{1})	
		 =\mbf{v}^{0} +\int_{a_{1}}^{s_{1}}\sum_{k=0}^{\delta_{1}-1}\bl(\mscr{g}_{1,k}^{\delta_{1}-1}\mbf{v}^{k+1}\br)(\theta)d\theta.
	\end{equation*}
	It follows that
	\begin{equation*}
		\tnf{b}^{-\delta_{1}}\partial_{s_{1}}^{0}\mbf{u}=\mbf{u}(a_{1})=\mbf{v}^{0},\qquad\text{and}\qquad
		\partial_{s_{1}}\mbf{u}=\sum_{k=0}^{\delta_{1}-1}\bl(\mscr{g}_{1,k}^{\delta_{1}-1}\mbf{v}^{k+1}\br).
	\end{equation*}
	By the induction hypothesis, then, for all $0\leq k\leq \tilde{\delta}_{1}=\delta_{1}-1$,
	\begin{equation*}
		\mbf{v}^{k+1}=\tnf{b}_{1}^{k-(\delta_{1}-1)}\partial_{s_{1}}^{k}(\partial_{s_{1}}\mbf{u})=\tnf{b}_{1}^{(k+1)-\delta_{1}}\partial_{s_{1}}^{k+1}\mbf{u},
	\end{equation*}
	whence $\mbf{v}^{k}=\tnf{b}_{1}^{k-\delta_{1}}\partial_{s_{1}}^{k}\mbf{u}$ for every $0\leq k\leq\delta_{1}$.
	Thus, for $\delta_{1}=\tilde{\delta}_{1}+1\in\N_{0}$, both implications of the lemma hold. By induction, the lemma holds for all $\delta_{1}\in\N_{0}$, and hence all $\delta\in\N_{0}^{N}$.
\end{pf}

\begin{example}
	Consider $\mbf{f}\in S_{2}^{(2,1)}[[a,b]\times[c,d]]$. Then
	\begin{align*}
		\tnf{b}_{1}^{2}\mbf{f}(y)&=\mbf{f}(a,y),	&
		\tnf{b}_{1}^{1}\partial_{x}\mbf{f}(y)&=\partial_{x}\mbf{f}(a,y),	&
		\tnf{b}_{1}^{0}\partial_{x}^2\mbf{f}(x,y)&=\partial_{x}^2\mbf{f}(x,y),	\\
		& &
		\tnf{b}_{2}^{1}\mbf{f}(x)&=\mbf{f}(x,c),	&
		\tnf{b}_{2}^{0}\partial_{y}\mbf{f}(x,y)&=\partial_{y}\mbf{f}(x,y).
	\end{align*}
	By Lemma~\ref{lem:fundamental_expansion_1D}, $\mbf{f}$ can be expanded as
	\begin{equation*}
		\mbf{f}(x,y)=\!\sum_{k=0}^{2}(\mscr{g}_{1,k}^{2}\tnf{b}_{1}^{k-2}\partial_{x}^{k}\mbf{f})(x,y)
		=\mbf{f}(a,y) +(x-a)\partial_{x}\mbf{f}(a,y) +\!\int_{a}^{x}\!(x-\theta)\partial_{x}^2\mbf{f}(\theta,y) d\theta,
	\end{equation*}
	or alternatively
	\begin{equation*}
		\mbf{f}(x,y)=\sum_{k=0}^{1}(\mscr{g}_{2,k}^{1}\tnf{b}_{2}^{k-1}\partial_{x}^{k}\mbf{f})(x,y)
		=\mbf{f}(x,c) +\!\int_{c}^{y}\!\partial_{y}\mbf{f}(x,\theta) d\theta.
	\end{equation*}
	Combining these expansions, we retrieve the expansion of $\mbf{u}$ in terms of all of its partial derivatives $D^{(i,j)}\mbf{u}=\partial_{x}^{i}\partial_{y}^{j}\mbf{u}$ for $(0,0)\leq (i,j)\leq (2,1)$ from Theorem~\ref{thm:fundamental_expansion_ND}, as
	\begin{align}\label{eq:expansions_example}
		&\mbf{u}(x,y) = \mbf{u}(a,c)
		+ (x\!-\!a)\partial_{x}\mbf{u}(a,c)   
		+\int_{a}^{x}\!(x\!-\!\theta)\partial_{x}^2\mbf{u}(\theta,c)d\theta	\\ 
		&\quad +\int_{c}^{y}\!\partial_{y}\mbf{u}(a,\eta)d\eta 
		+\!\int_{c}^{y}\!(x\!-\!a)\partial_{x}\partial_{y}\mbf{u}(a,\eta)d\eta     +\!\int_{a}^{x}\!\!\int_{c}^{y}(x\!-\!\theta)
		\partial_{x}^2\partial_{y}\mbf{u}(\theta,\eta)d\eta d\theta.	\notag
	\end{align}
\end{example}

Lemma~\ref{lem:fundamental_expansion_1D} shows that a function $\mbf{u}\in S_{2}^{\delta}[\Omega]$ can be expressed in terms of its partial derivatives $\partial_{s_{i}}^{\alpha_{i}}\mbf{u}$ with respect to a single variable $s_{i}$, using a suitable integral operator along $\Omega_{i}=[a_{i},b_{i}]$. 
To prove Theorem~\ref{thm:fundamental_expansion_ND}, now, this expansion is applied separately along each dimension of the hyperrectangle $\Omega:=\prod_{i=1}^{N}[a_{i},b_{i}]$, expanding $\mbf{u}\in S_{2}^{\delta}[\Omega]$ in terms of all of its admissible derivatives $D^{\alpha}\mbf{u}$ up to order $\alpha=\delta$.

\begin{pot}
	To prove that the theorem holds for arbitrary $\delta\in\N_{0}^{N}$, let $J\in\{0,\hdots,N\}$, and suppose that $\delta_{i}=0$ for all $i\in\{J+1,\hdots,N\}$. Then, it suffices to show that the theorem holds for $J=N$, which is done through induction on the value of $J$.

	\tbf{Base case }$\mbs{J=0:}$
	For $J=0$, we have $\delta=0^{N}$, so that $G_{\alpha}^{\delta}B^{\alpha-\delta}D^{\alpha}=I$ for all $0^{N}\leq \alpha\leq\delta=0^{N}$. In this case, the result holds trivially.
	
	\tbf{Induction step }$\mbs{J>1:}$
	Suppose for induction that the theorem holds for some $\tilde{J}\in\{0,\hdots,N-1\}$, and let $J=\tilde{J}+1$. Let $\delta\in\N_{0}^{N}$ be such that $\delta_{i}=0$ for $i> J$, and let $\tilde{\delta}=\delta-\delta_{J}\text{e}_{J}\in\N_{0}^{N}$, so that $\tilde{\delta}_{i}=\delta_{i}$ for $i\leq J-1$, and $\tilde{\delta}_{i}=0$ for $i> J-1$. Then, for every $0^{N}\leq \tilde{\alpha}\leq\tilde{\delta}$,
	\begin{equation*}
		\mcl{G}_{\tilde{\alpha}}^{\tilde{\delta}}=\prod_{i=1}^{J-1}\mscr{g}_{i,\tilde{\alpha}_{i}}^{\delta_{i}},\qquad
		B^{\tilde{\alpha}-\tilde{\delta}}=\prod_{i=1}^{J-1}\tnf{b}_{i}^{\tilde{\alpha}_{i}-\delta_{i}},\qquad
		D^{\tilde{\alpha}}=\prod_{i=1}^{J-1}\partial_{s_{i}}^{\tilde{\alpha}_{i}}.
	\end{equation*}
	Now, fix arbitrary $\mbf{u}\in S_{2}^{\delta}[\Omega]\subseteq S_{2}^{\tilde{\delta}}[\Omega]$. Invoking the induction hypothesis,
	and applying Lemma~\ref{lem:fundamental_expansion_1D}, it follows that
	\begin{align*}
		\mbf{u}&=\sum_{0^{N}\leq \tilde{\alpha}\leq \tilde{\delta}}\bl(\mcl{G}^{\tilde{\delta}}_{\tilde{\alpha}}B^{\tilde{\alpha}-\tilde{\delta}}D^{\tilde{\alpha}}\mbf{u}\br)	\\
		&=\sum_{0^{N}\leq \tilde{\alpha}\leq \tilde{\delta}}\left(\mcl{G}^{\tilde{\delta}}_{\tilde{\alpha}}B^{\tilde{\alpha}-\tilde{\delta}}D^{\tilde{\alpha}}\bbbl(\sum_{\alpha_{J}=0}^{\delta_{J}}\bl(\mscr{g}^{\delta_{J}}_{J,\alpha_{J}}\tnf{b}_{J}^{\alpha_{J}-\delta_{J}}\partial_{s_{J}}^{\alpha_{J}}\mbf{u}\br)\bbbr)\right)	\\
		&=\sum_{0^{N}\leq \alpha\leq \delta}\bbbbl(\prod_{i=1}^{J}\mscr{g}_{i,\alpha_{i}}^{\delta_{i}}\bbbbr)\bbbbl(\prod_{i=1}^{J}\tnf{b}_{i}^{\alpha_{i}-\delta_{i}}\bbbbr)\bbbbl(\prod_{i=1}^{J}\partial_{s_{i}}^{\alpha_{i}}\bbbbr)\mbf{u} 
		\hspace*{0.9cm}=\sum_{0^{N}\leq \alpha\leq \delta}\bl(\mcl{G}^{\delta}_{\alpha}B^{\alpha-\delta}D^{\alpha}\mbf{u}\br).
	\end{align*}
	Hence $\mbf{u}$ satisfies~\eqref{eq:fundamental_expansion_ND}. To see that the implication given by~\eqref{eq:inverse_expansion_ND} also holds, let $\{\mbf{v}^{\alpha}\in L_2[\Omega^{\alpha-\delta}]\mid 0^{N}\leq \alpha\leq \delta\}$  be such that
	\begin{align*}
		\mbf{u}&=\sum_{0^{N}\leq \alpha\leq \delta}\bl(\mcl{G}_{\alpha}^{\delta}\mbf{v}^{\alpha}\br)
		=\sum_{0^{N}\leq \tilde{\alpha}\leq\tilde{\delta}}\left(\mcl{G}_{\tilde{\alpha}}^{\tilde{\delta}}\bbbbl(\sum_{\alpha_{J}=0}^{\delta_{J}}\mscr{g}_{J,\alpha_{J}}^{\delta_{J}}\mbf{v}^{\alpha}\bbbbr)\right),
	\end{align*}
	where $\alpha=\tilde{\alpha}+\alpha_{J}\text{e}_{J}$.
	Then, by the induction hypothesis,
	\begin{equation*}
		\sum_{\alpha_{J}=0}^{\delta_{J}}\mscr{g}_{J,\alpha_{J}}^{\delta_{J}}\mbf{v}^{\alpha}
		=B^{\tilde{\alpha}-\tilde{\delta}}D^{\tilde{\alpha}}\mbf{u},\qquad 0^{N}\leq \tilde{\alpha}\leq \tilde{\delta}.
	\end{equation*}
	By Lemma~\ref{lem:fundamental_expansion_1D}, it follows that for all $0^{N}\leq \alpha\leq\delta$,
	\begin{equation*}
		\mbf{v}^{\alpha}
		=\tnf{b}_{J}^{\alpha_{J}-\delta_{J}}\partial_{s_{J}}^{\alpha_{J}}\bbl(B^{\tilde{\alpha}-\tilde{\delta}}D^{\tilde{\alpha}}\mbf{u}\bbr)
		=B^{\alpha-\delta}D^{\alpha}\mbf{u}.
	\end{equation*}
	Thus, for the given $J=\tilde{J}+1\in\{0,\hdots,N\}$, both implications of the theorem hold. By induction, the theorem holds for $J=N$, and hence for all $\delta\in\N_{0}^{N}$.
\end{pot}

Theorem~\ref{thm:fundamental_expansion_ND} shows that a differentiable function $\mbf{u}$ is uniquely defined by its derivatives $D^{\alpha}\mbf{u}$ for $0^{N}\leq\alpha\leq\delta$, evaluated at the lower boundaries $\Omega^{\alpha-\delta}$ of the domain. In particular, a derivative $D^{\alpha}\mbf{u}$ in this expansion is evaluated at the boundary $s_{i}=a_{i}$ for some $i\in\{1,\hdots,N\}$ if and only if $\alpha_{i}\neq\delta_{i}$. As a consequence, for each dimension $i\in\{1,\hdots,N\}$, the boundary value $B^{\alpha-\delta}D^{\alpha}\mbf{u}$ is a function of $s_{i}\in[a_{i},b_{i}]$ if and only if $\alpha_{i}=\delta_{i}$ -- i.e. if the highest-order derivative is taken with respect to $s_{i}$.
Therefore, the boundary values $B^{\alpha-\delta}D^{\alpha}\mbf{u}\in L_{2}[\Omega^{\alpha-\delta}]$ in the expansion of $\mbf{u}\in S^{\delta}_{2}[\Omega]$ need not be differentiable, and are free of the regularity constraints imposed upon $\mbf{u}$. An associated aggregate space for these boundary values can be defined as
\begin{equation*}
	\mbf{L}_{2}^{\delta}[\Omega]:=\prod_{0^{N}\leq \alpha\leq \delta} L_{2}[\Omega^{\alpha-\delta}],
\end{equation*}
with norm $\|\mbf{v}\|_{\mbf{L}_{2}^{\delta}}^2:=\sum_{0^{N}\leq\alpha\leq\delta} \|\mbf{v}^{\alpha}\|_{L_2}^2$ for $\mbf{v}=(\mbf{v}^{\alpha})\in \mbf{L}_{2}^{\delta}[\Omega]$, with $\mbf{v}^{\alpha}\in L_{2}[\Omega^{\alpha-\delta}]$ for $0^{N}\leq\alpha\leq\delta$. For any $\mbf{u}\in S_{2}^{\delta}[\Omega]$, then, there exists an associated $(\mbf{v}^{\alpha})\in \mbf{L}_{2}^{\delta}[\Omega]$, defined by $\mbf{v}^{\alpha}=B^{\alpha-\delta}D^{\alpha}\mbf{u}$ for each $\alpha$. Conversely, by Theorem~\ref{thm:fundamental_expansion_ND}, for any $(\mbf{v}^{\alpha})\in \mbf{L}_{2}^{\delta}[\Omega]$ there exists a unique $\mbf{u}=\sum_{0^{N}\leq\alpha\leq\delta}\mcl{G}^{\delta}_{\alpha}\mbf{v}^{\alpha}\in S_{2}^{\delta}[\Omega]$ such that $\mbf{v}^{\alpha}=B^{\alpha-\delta}D^{\alpha}\mbf{u}$ for all $\alpha$. As such, there exists a bijective map between $\mbf{L}_{2}^{\delta}[\Omega]$ and $S_{2}^{\delta}[\Omega]$, as shown in the following corollary.

\begin{cor}\label{cor:bijection}
	Let $\delta\in\N_{0}^{N}$ and $\Omega=\prod_{i=1}^{N}[a_{i},b_{i}]\subseteq\R^{N}$. Define the operator $\mcl{G}^{\delta}:\mbf{L}_{2}^{\delta}[\Omega]\to S^{\delta}_{2}[\Omega]$ by
	\begin{equation*}
		\quad\mcl{G}^{\delta}\mbf{v}:=\sum_{0^{N}\leq\alpha\leq\delta}\mcl{G}^{\delta}_{\alpha}\mbf{v}^{\alpha},\qquad
		\forall \mbf{v}=(\mbf{v}^{\alpha})\in\mbf{L}_{2}^{\delta}[\Omega],
	\end{equation*}
	where $\mbf{v}^{\alpha}\in L_2[\Omega^{\alpha-\delta}]$ and $\mcl{G}^{\delta}_{\alpha}$ is as defined in Theorem~\ref{thm:fundamental_expansion_ND} for all $0^{N}\leq\alpha\leq\delta$. Then, $\mcl{G}^{\delta}:\mbf{L}_{2}^{\delta}[\Omega]\to S^{\delta}_{2}[\Omega]$ is bounded and bijective.
\end{cor}
\begin{pf}
	We first show that $\mcl{G}^{\delta}:\mbf{L}_{2}^{\delta}[\Omega]\to S^{\delta}_{2}[\Omega]$ is surjective. Fix arbitrary $\mbf{u}\in S^{\delta}_{2}[\Omega]$, and let $\mbf{v}^{\alpha}:=B^{\alpha-\delta}D^{\alpha}\mbf{u}\in L_2[\Omega^{\alpha-\delta}]$ for all $0^{N}\leq\alpha\leq\delta$. Then, $\mbf{v}=(\mbf{v}^{\alpha})\in \mbf{L}_{2}^{\delta}[\Omega]$, and so by Theorem~\ref{thm:fundamental_expansion_ND}
	\begin{equation*}
		\mbf{u} =\sum_{0^{N}\leq\alpha\leq\delta}\mcl{G}^{\delta}_{\alpha} B^{\alpha-\delta}D^{\alpha}\mbf{u}
		=\sum_{0^{N}\leq\alpha\leq\delta}\mcl{G}^{\delta}_{\alpha} \mbf{v}^{\alpha}
		=\mcl{G}^{\delta}\mbf{v}.
	\end{equation*}
	It follows that the map $\mcl{G}^{\delta}:\mbf{L}_{2}^{\delta}[\Omega]\to S^{\delta}_{2}[\Omega]$ is surjective.
	
	Next, to prove that $\mcl{G}^{\delta}$ is injective, suppose that $\mcl{G}^{\delta}\mbf{v}=\mbf{0}\in S^{\delta}_{2}[\Omega]$ for some $\mbf{v}=(\mbf{v}^{\alpha})\in \mbf{L}_{2}^{\delta}[\Omega]$ with $\mbf{v}^{\alpha}\in L_{2}[\Omega^{\alpha-\delta}]$ for $0^{N}\leq\alpha\leq\delta$. By Theorem~\ref{thm:fundamental_expansion_ND}, then, $\mbf{v}^{\alpha}=B^{\alpha-\delta}D^{\alpha}\mbf{0}=\mbf{0}$ for every $0^{N}\leq\alpha\leq\delta$, and thus $\mbf{v}=\mbf{0}\in \mbf{L}_{2}^{\delta}[\Omega]$. It follows that $\tnf{ker}(\mcl{G}^{\delta})=\{\mbf{0}\}$, and therefore $\mcl{G}^{\delta}:\mbf{L}_{2}^{\delta}[\Omega]\to S^{\delta}_{2}[\Omega]$ is injective.
	
	Finally, to show that $\mcl{G}^{\delta}:\mbf{L}_{2}^{\delta}[\Omega]\to S^{\delta}_{2}[\Omega]$ is bounded, fix again arbitrary $\mbf{v}=(\mbf{v}^{\alpha})\in \mbf{L}_{2}^{\delta}[\Omega]$, and let $\mbf{u}=\mcl{G}^{\delta}\mbf{v}\in S^{\delta}_{2}[\Omega]$. Then $\mbf{v}^{\alpha}=B^{\alpha-\delta}D^{\alpha}\mbf{u}\in L_2[\Omega^{\alpha-\delta}]$ for every $0^{N}\leq\alpha\leq\delta$. 
	Since $D^{\gamma}\mbf{u}\in S_{2}^{\delta-\gamma}[\Omega]$ for every $0^{N}\leq\gamma\leq\delta$, by Theorem~\ref{thm:fundamental_expansion_ND}, $D^{\gamma}$ can be expanded as 
	\begin{equation*}
		D^{\gamma}\mbf{u}=\!\!\sum_{0^{N}\leq\alpha\leq\delta-\gamma}\!\!\!G^{\delta-\gamma}_{\alpha}B^{\alpha-(\delta-\gamma)}D^{\alpha+\gamma}\mbf{u}
		=\sum_{\gamma\leq\alpha\leq\delta}G^{\delta-\gamma}_{\alpha-\gamma} B^{\alpha-\delta}D^{\alpha}\mbf{u}
		=\sum_{\gamma\leq\alpha\leq\delta}G^{\delta-\gamma}_{\alpha-\gamma}\mbf{v}^{\alpha}.
	\end{equation*}
	Now, since each of the operators $\mcl{G}^{\delta-\gamma}_{\alpha-\gamma}$ is either a multiplier operator defined by a bounded (polynomial) function, or an integral operator defined by a square-integrable (polynomial) kernel, the operators $\mcl{G}^{\delta-\gamma}_{\alpha-\gamma}:L_{2}[\Omega^{\alpha-\delta}]\to L_{2}[\Omega]$ are bounded for every $0^{N}\leq\gamma\leq\alpha\leq\delta\in \N_{0}^{N}$. It follows that
	\begin{align*}
		\|D^{\gamma}\mbf{u}\|_{L_{2}}
		&=\norm{\sum_{\gamma\leq\alpha\leq\delta}G^{\delta-\gamma}_{\alpha-\gamma}\mbf{v}^{\alpha}}_{L_2}	\\
		&\leq \sum_{\gamma\leq\alpha\leq\delta}\|G^{\delta-\gamma}_{\alpha-\gamma}\mbf{v}^{\alpha}\|_{L_2}	\\
		&\leq \sum_{\gamma\leq\alpha\leq\delta}\|G^{\delta-\gamma}_{\alpha-\gamma}\|_{\tnf{op}}\|\mbf{v}^{\alpha}\|_{L_2} 	
		\leq C_{\delta,\gamma}\sum_{0^{N}\leq\alpha\leq\delta}\|\mbf{v}^{\alpha}\|_{L_2}
		\leq C_{\delta,\gamma}\sqrt{\rho_{\delta}}\|\mbf{v}\|_{\mbf{L}_{2}},
	\end{align*}
	where we define $\rho_{\delta}:=\prod_{i=1}^{N}(\delta_{i}+1)$ and $C_{\delta,\gamma}:=\max_{\gamma\leq\alpha\leq\delta}\|G^{\delta-\gamma}_{\alpha-\gamma}\|_{\tnf{op}}$, for every $0^{N}\leq\gamma\leq\delta$. Finally,
	\begin{equation*}
		\|\mbf{u}\|_{S_{2}^{\delta}}^2
		:=\sum_{0^{N}\leq\gamma\leq\delta}\|D^{\gamma}\mbf{u}\|_{L_{2}}^2
		\leq \sum_{0^{N}\leq\gamma\leq\delta}\rho_{\delta} C_{\delta,\gamma}^2\|\mbf{v}\|_{\mbf{L}_{2}}^2
		= \rho_{\delta}C_{\delta}^2 \|\mbf{v}\|_{\mbf{L}_{2}}^2,
	\end{equation*}
	where we define $C_{\delta}^2 := \sum_{0^{N}\leq\gamma\leq\delta} C_{\delta,\gamma}^2$.
	Thus $\mcl{G}^{\delta}:\mbf{L}_{2}^{\delta}[\Omega]\to S^{\delta}_{2}[\Omega]$ is bounded, with $\|\mcl{G}^{\delta}\|_{\tnf{op}}\leq \sqrt{\rho_{\delta}}C_{\delta}$.
\end{pf}

\begin{example}
	Consider again $\mbf{f}\in S_{2}^{(2,1)}[[a,b]\times[c,d]]$. Then
	\begin{equation*}
		\mbf{L}_{2}^{(2,1)}[[a,b]\times[c,d]]
		=\R^{2}\times L_2[[a,b]]\times L_2^{2}[[c,d]]\times L_2[[a,b]\times [c,d]].
	\end{equation*}
	By Corollary~\ref{cor:bijection}, there exists a unique element $\mbf{h}=(\mbf{h}^{\alpha})\in \mbf{L}_{2}^{(2,1)}[[a,b]\times[c,d]]$ such that $\mbf{f}=\mcl{G}^{(2,1)}\mbf{h}$. In particular, by the proof of the corollary, $\mbf{h}^{\alpha}=B^{\alpha-(2,1)}D^{\alpha}\mbf{f}$ for every $(0,0)\leq\alpha\leq (2,1)$, so that
	\begin{align*}
		&\mbf{h}(x,y)
		=\bmat{\mbf{h}^{(0,0)}&\mbf{h}^{(1,0)}&\mbf{h}^{(2,0)}(x)&\mbf{h}^{(0,1)}(y)&\mbf{h}^{(1,1)}(y)&\mbf{h}^{(2,1)}(x,y)}^T	\\
		&=\bmat{\mbf{f}(a,c)&\partial_{x}\mbf{f}(a,c)&\partial_{x}^{2}\mbf{f}(x,c)&\partial_{y}\mbf{f}(a,y)&\partial_{x}\partial_{y}\mbf{f}(a,y)&\partial_{x}^{2}\partial_{y}\mbf{f}(x,y)}^T.
	\end{align*}
	Then, $\mbf{f}$ can be reconstructed from $\mbf{h}$ as
	\begin{align*}
			&\mbf{f}(x,y)=(\mcl{G}^{(2,1)}\mbf{h})(x,y)	\\
			&=\mbf{h}^{(0,0)} +(x-a)\mbf{h}^{(1,0)} +\int_{a}^{x}(x-\theta)\mbf{h}^{(2,0)}(\theta)d\theta \\
			&\quad+\int_{c}^{y}\mbf{h}^{(0,1)}(\eta)d\eta +(x-a)\int_{c}^{y}\mbf{h}^{(1,1)}(\eta)d\eta +\int_{a}^{x}\!\int_{c}^{y}(x-\theta)\mbf{h}^{(2,1)}(\theta,\eta)d\eta d\theta.
	\end{align*}
\end{example}

Corollary~\ref{cor:bijection} proves that any differentiable function $\mbf{u}\in S_{2}^{\delta}[\Omega]$ is uniquely defined by an associated element $\mbf{v}\in \mbf{L}_{2}^{\delta}[\Omega]$, free of the continuity constraints imposed upon $\mbf{u}$. In the next section, this relation will be used to perform polynomial approximation of functions $\mbf{u}\in S_{2}^{\delta}[\Omega]$ in the continuity-constraints-free space $\mbf{L}_{2}^{\delta}[\Omega]$.

\begin{remark}
	Although Theorem~\ref{thm:fundamental_expansion_ND} offers an expansion of $\mbf{u}\in S^{\delta}_{2}[\Omega]$ in terms of its derivatives $D^{\alpha}\mbf{u}$ evaluated at the lower boundaries $\Omega^{\alpha-\delta}$ of the domain, a similar expansion may readily be given in terms of the derivatives $D^{\alpha}\mbf{u}$ along upper boundaries $\Omega^{\delta-\alpha}$, for $0^{N}\leq\alpha\leq\delta$.
 
For example, just as $\mbf{u}\in S_{2}^{(2,1)}[[a,b]\times[c,d]]$ can be expressed in terms of its derivatives at $x=a$ and $y=c$ as in~\eqref{eq:expansions_example}, so too can it be expressed in terms of its derivatives at $x=b$ and $y=d$ as
\begin{align*}
	&\mbf{u}(x,y) = \mbf{u}(b,d)
	- (b\!-\!x)\partial_{x}\mbf{u}(b,d)  
	-\int_{x}^{b}\!\!(x\!-\!\theta)\partial_{x}^2\mbf{u}(\theta,d)d\theta	\\ 
	&\quad-\int_{y}^{d}\!\!\partial_{y}\mbf{u}(b,\eta)d\eta
	+\!\int_{y}^{d}\!\!(b\!-\!x)\partial_{x}\partial_{y}\mbf{u}(b,\eta)d\eta     +\!\int_{x}^{b}\!\!\int_{y}^{d}\!(x\!-\!\theta)
	\partial_{x}^2\partial_{y}\mbf{u}(\theta,\eta)d\eta d\theta.
\end{align*}
Such different expansions of $\mbf{u}\in S_{2}^{\delta}$ in terms of different boundary values are of particular interest when considering functions constrained by boundary conditions, in which case the value of certain terms $B^{\beta}D^{\alpha}\mbf{u}$ is known. Given sufficient and suitable boundary conditions, then, Theorem~\ref{thm:fundamental_expansion_ND} proves that there exists a unique map between $\mbf{u}$ and its highest-order derivative $D^{\delta}\mbf{u}$. A full derivation of such a relation is left for future works.
\end{remark}

\section{Polynomial Approximation of Functions in $S_{2}^{\delta}$}\label{sec:approximation}

Theorem~\ref{thm:fundamental_expansion_ND} offers an alternative representation of functions $\mbf{u}\in S_{2}^{\delta}[\Omega]$ on the hyperrectangle $\Omega=\prod_{i=1}^{N}[a_{i},b_{i}]$ in terms of their derivatives $D^{\alpha}\mbf{u}$ along lower boundaries $\Omega^{\alpha-\delta}$ of the domain. 
The benefit of this representation is that the boundary values $(B^{\alpha-\delta}D^{\alpha}\mbf{u})^{0^{N}\leq\alpha\leq\delta}\in \mbf{L}_{2}^{\delta}[\Omega]$ are all elements of $L_{2}$ rather than $S_{2}^{\delta}$, and therefore not required to be differentiable or even continuous. 
Using the bijection $\mcl{G}^{\delta}:\mbf{L}_{2}^{\delta}[\Omega]\to S^{\delta}_{2}[\Omega]$ defined in Corollary~\ref{cor:bijection}, then, analysis, approximation, and simulation of functions in $S_{2}^{\delta}[\Omega]$ may be performed in the less restrictive space $\mbf{L}_{2}^{\delta}[\Omega]$.

In this section, we illustrate one application of Theorem~\ref{thm:fundamental_expansion_ND}, concerning polynomial approximation of functions in $S_{2}^{\delta}[\Omega]$. In particular, 
approximation of $\mbf{u}\in S_{2}^{\delta}$ can be performed by first projecting each of the derivatives $B^{\alpha-\delta}D^{\alpha}\mbf{u}$ onto a polynomial subspace, and subsequently using Theorem~\ref{thm:fundamental_expansion_ND} to reconstruct an approximation of $\mbf{u}$.
The following subsection shows how this may be achieved using a basis of Legendre polynomials, to obtain polynomial approximations of functions in $S_{2}^{\delta}$. In the subsequent subsection, this approach is then adapted to use a basis of step functions, to construct piecewise polynomial approximations. Finally, convergence of approximations using both methods is demonstrated for 1D and 2D examples in Subsection~\ref{sec:approximation:examples}.

Throughout this section, the domain $\Omega=\prod_{i=1}^{N}[-1,1]$ is assumed to be a hypercube for ease of notation.

\subsection{Approximation using Legendre Polynomials}\label{sec:approximation:legendre}

As a first method for approximation of functions $\mbf{u}\in S_{2}^{\delta}[\Omega]$, 
consider the projection of such functions onto the space of polynomials. In particular, for $d\in\N_{0}^{N}$, let $\R^{d}[\Omega]$ denote the space of real-valued polynomials of degree at most $d_{i}$ in each variable $s_{i}\in[-1,1]$ for $i\in\{1,\hdots,N\}$. A polynomial approximation $P_{d}\mbf{u}\in\R^{d}[\Omega]$ of a function $\mbf{u}\in S_{2}^{\delta}[\Omega]$ can be computed using Legendre polynomials, which may be recursively defined for $x\in[-1,1]$ as
\begin{align*}
	&\varphi^{0}(x)=1,\hspace*{3.0cm} \varphi^{1}(x)=x,	\\
	&\varphi^{k+1}(x)=\frac{2k+1}{k+1}\varphi^{k}(x)-\frac{k}{k+1}\varphi^{k-1}(x),\quad k\geq 1,
\end{align*}
and which are orthogonal with respect to the standard $L_{2}$ inner product, satisfying
\begin{equation*}
	\int_{-1}^{1}\varphi^{m}(x)\varphi^{n}(x)dx
	=\begin{cases}
		\frac{1}{m+\frac{1}{2}},	&	m=n,	\\
		0,		&	\text{else}.
	\end{cases}
\end{equation*}
Letting $\mbs{\phi}^{d}(s):=\prod_{i=1}^{N}\sqrt{d_{i}+\frac{1}{2}}\:\varphi^{d_{i}}(s_{i})$ for $d\in\N_{0}^{N}$ and $s\in\Omega$, a function $\mbf{u}\in S_{2}^{\delta}[\Omega]$ can then be projected onto the space $\R^{d}[\Omega]$ as 
\begin{equation}\label{eq:projection_leg}
	P_{d}\mbf{u}=\!\sum_{0^{N}\leq\alpha\leq d}\! c_{\alpha}\mbs{\phi}^{\alpha},\qquad \text{where}\qquad
	c_{\alpha}:=\int_{\Omega}\mbf{u}(s)\mbs{\phi}^{\alpha}(s)ds.
\end{equation}
Using this approach, it is well-known the projection $P_{d}\mbf{u}$ offers an $L_{2}$-optimal polynomial approximation of $\mbf{u}$, in the sense that
\begin{equation}\label{eq:proj_optimal_L2}
	\|\mbf{u}-P_{d}\mbf{u}\|_{L_2}=\min_{\mbf{v}\in\R^{d}[\Omega]}\|\mbf{u}-\mbf{v}\|_{L_2}.
\end{equation}
Unfortunately, however, this projection $P_{d}\mbf{u}$ may not accurately approximate the function $\mbf{u}$ in terms of the Sobolev norm. In particular, suppose that $\delta=\delta_{0}\cdot 1^{N}$ for some $\delta_{0}\in\N$, and let $\mbf{u}\in S_{2}^{\delta}[\Omega]$. Define the standard Sobolev norm of $\mbf{u}$ as
\begin{equation*}
	\|\mbf{u}\|_{W_{2}^{\delta_{0}}}^2:=\sum_{\alpha\in\N_{0}^{N},~\|\alpha\|_{1}\leq\delta_{0}}\|D^{\alpha}\mbf{u}\|_{L_{2}}^2.
\end{equation*}
Clearly then, $\|\mbf{u}\|_{W_{2}^{\delta_{0}}}\leq\|\mbf{u}\|_{S_{2}^{\delta}}$, and therefore $\mbf{u}\in W_{2}^{\delta_{0}}[\Omega]$. Under these conditions, it has been shown that if $d=d_{0}\cdot 1^{N}$ for some $d_{0}\in\N$, then the $L_{2}$ norm of the error in the projection $P_{d}\mbf{u}$ is bounded as~\cite{canuto1982approximation}
\begin{equation*}
	\|\mbf{u}-P_{d}\mbf{u}\|_{L_2}\leq C_{0} d_{0}^{-\delta_{0}} \|\mbf{u}\|_{W_{2}^{\delta_{0}}} \leq C_{0} d_{0}^{-\delta_{0}}\|\mbf{u}\|_{S_{2}^{\delta}},
\end{equation*}
for some constant $C_{0}$, proving that the error decays as at least $\mcl{O}(d_{0}^{-\delta_{0}})$ with the degree $d_{0}$ of the polynomial projection.
However, the Sobolev norm of the error is bounded only as~\cite{canuto2007projection,gottlieb1977spectral_methods}
\begin{equation*}
	\|\mbf{u}-P_{d}\mbf{u}\|_{W_2^{\gamma_{0}}}\leq C_{\gamma_{0}} d_{0}^{2\gamma_{0}-\delta_{0}-\frac{1}{2}} \|\mbf{u}\|_{W_{2}^{\delta_{0}}},\qquad 0\leq \gamma_{0}\leq \delta_{0}\in\N,
\end{equation*}
for some constant $C_{\gamma_{0}}$, not guaranteeing any convergence of the approximation in the Sobolev norm.
Indeed, the Sobolev norm $\|\mbf{u}-P_{d}\mbf{u}\|_{W_2^{\delta_{0}}}$ of the error need only be bounded as $\mcl{O}(d_{0}^{\delta_{0}-1/2})$, and this error may actually increase with the degree $d_{0}$ of the polynomials -- see e.g. Fig.~\ref{fig:invsqrt_conv_Legendre} in Subsection~\ref{sec:approximation:examples}.

The problem, of course, is that using the $L_{2}$ inner product to project $\mbf{u}$ onto a polynomial subspace, the resulting projection $P_{d}\mbf{u}$ need not capture any of the derivatives $D^{\alpha}\mbf{u}$ of the function.
Instead, accurate approximation in terms of a Sobolev norm requires projecting not just the function $\mbf{u}$, but also its derivatives $D^{\alpha}\mbf{u}$. Although this can be achieved using a Sobolev inner product with suitable polynomial basis~\cite{marcellan1993sobolev_polynomials,marcellan2015sobolev_polynomials,xu2018sobolev_polynomials}, we propose to instead use Theorem~\ref{thm:fundamental_expansion_ND}, to reconstruct an approximation of $\mbf{u}$ from $L_{2}$ projections of its derivatives $\mbf{v}^{\alpha}=B^{\alpha-\delta}D^{\alpha}\mbf{u}$. In particular, if $\mbf{u}=\sum_{0^{N}\leq\alpha\leq\delta}\mcl{G}_{\alpha}^{\delta}\mbf{v}^{\alpha}$, each function $\mbf{v}^{\alpha}\in L_2[\Omega^{\alpha-\delta}]$ can be projected onto $\R^{d}[\Omega^{\alpha-\delta}]$ as
\begin{equation*}
	P_{\kappa_{d}^{\alpha-\delta}}\mbf{v}^{\alpha}=\!\!\sum_{0^{N}\leq\gamma\leq \kappa_{d}^{\alpha-\delta}}\!\!\! c_{\gamma}\mbs{\phi}^{\gamma},\qquad \text{where}\qquad
	c_{\gamma}:=\int_{\Omega^{\alpha-\delta}}\!\!\mbf{u}(s)\mbs{\phi}^{\gamma}(s)ds,
\end{equation*}
where for $-\delta\leq\beta\leq 0^{N}$ we define $\kappa_{d}^{\beta}\in\N_{0}^{N}$ by
\begin{equation*}
	(\kappa_{d}^{\beta})_{i}:=
	\begin{cases}
		d_{i},	&	\beta_{i}=0,	\\
		0,		&	\text{else},
	\end{cases}
	\hspace*{2.0cm} \forall i\in\{1,\hdots,N\},
\end{equation*}%
so that e.g. $P_{\kappa_{d}^{\alpha-\delta}}\mbf{v}^{\alpha}=c_{0^{N}}\mbs{\phi}^{0^{N}}=\mbf{v}^{\alpha}$ when $\mbf{v}^{\alpha}\in\R$ is finite-dimensional.
Then, for any $0^{N}\leq \gamma\leq\delta$ and $d\in\N_{0}^{N}$, an \textit{order-$\gamma$ Sobolev projection} of $\mbf{u}$ onto $\R^{d+\gamma}[\Omega]$ can be constructed as
\begin{equation}\label{eq:projection_leg_sob}
	P_{d}^{\gamma}\mbf{u}:=\sum_{0^{N}\leq\alpha\leq\gamma}\mcl{G}_{\alpha}^{\gamma}\:\bbl(P_{\kappa_{d}^{\alpha-\gamma}}\bl[B^{\alpha-\gamma}D^{\alpha}\mbf{u}\br]\bbr).
\end{equation}
Approximating $\mbf{u}\in S_{2}^{\delta}[\Omega]$ in this manner is equivalent to projecting the function onto a polynomial subspace using a discrete-continuous Sobolev inner product~\cite{diaz2021discrete_continuous_sobolev}. In particular, defining a discrete-continuous norm on
$S_{2}^{\delta}[\Omega]$ as the norm of its associated boundary elements in $\mbf{L}_{2}^{\delta}[\Omega]$,
\begin{equation}\label{eq:norm_sob_dc}
	\norm{\mbf{u}}_{\mbf{L}_{2}^{\delta}}^2
	:=\!\!\sum_{0^{N}\leq\alpha\leq\delta}\norm{B^{\alpha-\delta}D^{\alpha}\mbf{u}}_{L_{2}}^2,\qquad 
	\mbf{u}\in S_{2}^{\delta}[\Omega],
\end{equation}
the following proposition shows that the approximation $P_{d}^{\delta}\mbf{u}$ is optimal with respect to this norm.

\begin{prop}
	Let $\mbf{u}\in S_{2}^{\delta}[\Omega]$ for some $\delta\in\N_{0}^{N}$, and define the polynomial approximation $P_{d}^{\gamma}\mbf{u}$ for $0^{N}\leq\gamma\leq\delta$ as in~\eqref{eq:projection_leg_sob}. Then, for every $d\in\N_{0}^{N}$, 
	\begin{equation*}
		\|\mbf{u}-P_{d}^{\gamma}\mbf{u}\|_{\mbf{L}_{2}^{\gamma}}=\min_{\mbf{v}\in\R^{d+\gamma}[\Omega]}\|\mbf{u}-\mbf{v}\|_{\mbf{L}_{2}^{\gamma}},
	\end{equation*}
	where the norm $\|.\|_{\mbf{L}_{2}^{\gamma}}$ is as defined in~\eqref{eq:norm_sob_dc}.
\end{prop}

\begin{pf}
	Fix arbitrary $\mbf{u}\in S_{2}^{\delta}[\Omega]$ and $\mbf{v}\in\R^{d+\gamma}[\Omega]$. Then $B^{\alpha-\gamma}D^{\alpha}\mbf{u}\in L_{2}[\Omega^{\alpha-\gamma}]$ and $B^{\alpha-\gamma}D^{\alpha}\mbf{v}\in\R^{\kappa_{d}^{\alpha-\gamma}}[\Omega^{\alpha-\gamma}]$ for every $0^{N}\leq\alpha\leq\gamma$. 
	By optimality of the projection $P_{d}$ in $L_2[\Omega]$ (Eqn.~\eqref{eq:proj_optimal_L2}), it follows that
	\begin{align*}
		\|\mbf{u}-\mbf{v}\|_{\mbf{L}_{2}^{\gamma}}^2
		&=\!\sum_{0^{N}\leq\alpha\leq\gamma}\norm{B^{\alpha-\gamma}D^{\alpha}\mbf{u}-B^{\alpha-\gamma}D^{\alpha}\mbf{v}}_{L_{2}}^2	\\
		&\geq\! \sum_{0^{N}\leq\alpha\leq\gamma}\norm{B^{\alpha-\gamma}D^{\alpha}\mbf{u}-P_{\kappa_{d}^{\alpha-\gamma}}[B^{\alpha-\gamma}D^{\alpha}\mbf{u}]}_{L_{2}}^2	\\
		&=\!\sum_{0^{N}\leq\alpha\leq\gamma}\norm{B^{\alpha-\gamma}D^{\alpha}\bl(\mbf{u}-P_{d}^{\gamma}\mbf{u}\br)}_{L_{2}}^2 \\
		&=\norm{\mbf{u}-P_{d}^{\gamma}\mbf{u}}_{\mbf{L}_{2}^{\gamma}}^2,
	\end{align*}
	where $B^{\alpha-\gamma}D^{\alpha}(P_{d}^{\gamma}\mbf{u})=P_{\kappa_{d}^{\alpha-\gamma}}[B^{\alpha-\gamma}D^{\alpha}\mbf{u}]$ by Theorem~\ref{thm:fundamental_expansion_ND} and Eqn.~\eqref{eq:projection_leg_sob}.
\end{pf}

\subsection{Expansion using Step Function Basis}\label{sec:approximation:stepfunction}

In the previous subsection, it was shown how Theorem~\ref{thm:fundamental_expansion_ND} can be used to construct a polynomial approximation of a differentiable function $\mbf{u}\in S_{2}^{\delta}[\Omega]$, using polynomial projections of the boundary values $\mbf{v}^{\alpha}=B^{\alpha-\delta}D^{\alpha}\mbf{u}$.
However, computing these projections requires integrating the product $\mbf{v}^{\alpha}\mbs{\phi}^{\gamma}$ for each polynomial $\mbs{\phi}^{\gamma}$ up to degree $\gamma=d$, which becomes numerically challenging for large degrees $d\in\N_{0}^{N}$. Moreover, since the boundary values $\mbf{v}^{\alpha}\in L_2[\Omega^{\alpha-\delta}]$ -- and in particular the highest-order derivative $\mbf{v}^{\delta}=D^{\delta}\mbf{u}\in L_2[\Omega]$ -- need not be continuous, a polynomial projection of these functions may not be well-behaved, exhibiting e.g. Gibbs phenomenon.

As an alternative to the polynomial projection method presented in the previous subsection, it is also possible to perform projection using a basis of step functions. In particular, for $K\in\N^{N}$, decompose the hypercube $\Omega=[-1,1]^{N}$ into $\prod_{i=1}^{N} K_{i}$ disjoint hypercubes $\Gamma_{K}^{\lambda}\subseteq \Omega$ for $1^{N}\leq \lambda\leq K$ of equal dimensions, each with volume $\Delta_{K}=\prod_{i=1}^{N}\frac{2}{K_{i}}$, so that $\Omega=\bigcup_{1^{N}\leq \lambda\leq K}\Gamma_{K}^{\lambda}$. 
Then, we define a piecewise constant approximation $Q_{K}\mbf{v}$ of $\mbf{v}\in L_2[\Omega]$ as the projection of $\mbf{v}$ onto a basis of step functions on these cells, as
\begin{equation*}
	Q_{K}\mbf{v}\br|_{\Gamma_{K}^{\lambda}}
	=\frac{1}{\Delta_{K}}\!\int_{\Gamma_{K}^{\lambda}}\!\!\mbf{v}(s) ds,\qquad 1^{N}\leq\lambda\leq K.
\end{equation*}
Projecting each derivative $\mbf{v}^{\alpha}=B^{\alpha-\gamma}D^{\alpha}\mbf{u}\in L_2[\Omega^{\alpha-\gamma}]$ in this manner, a piecewise polynomial approximation of $\mbf{u}$ can then be reconstructed as
\begin{equation}\label{eq:projection_stepfun_sob}
	Q_{K}^{\gamma}\mbf{u}=\sum_{0^{N}\leq\alpha\leq\gamma}\mcl{G}_{\alpha}^{\gamma}Q_{K}\bl[B^{\alpha-\gamma}D^{\alpha}\mbf{u}\br]. \\[-0.6em]
\end{equation}
Note that, by expanding the derivative $D^{\delta}\mbf{u}\in L_2[\Omega]$ using a basis of step functions, the proposed projection $Q_{K}^{\delta}\mbf{u}\in S_{2}^{\delta}[\Omega]$ is sufficiently smooth to capture all of the derivatives of the function $\mbf{u}\in S_{2}^{\delta}[\Omega]$ -- unlike e.g. a direct step function expansion $Q_{K}^{0}\mbf{u}:=Q_{K}\mbf{u}\in L_2[\Omega]$ -- but will not suffer from Gibbs phenomenon -- unlike e.g. a polynomial projection. In addition, the integrals $\int_{\Gamma_{K}^{\lambda}}\mbf{v}^{\alpha}(s) ds$ in the step function projection are in general much easier to compute than the integrals $\int_{\Omega}\mbs{\phi}^{d}(s)\mbf{v}(s)ds$ for more complicated basis functions $\mbs{\phi}^{d}$, such as the Legendre polynomials.

\subsection{Numerical Examples}\label{sec:approximation:examples}

\subsubsection{1D Numerical Example}

\begin{figure}
	\centering
	\includegraphics[width=1.0\textwidth]{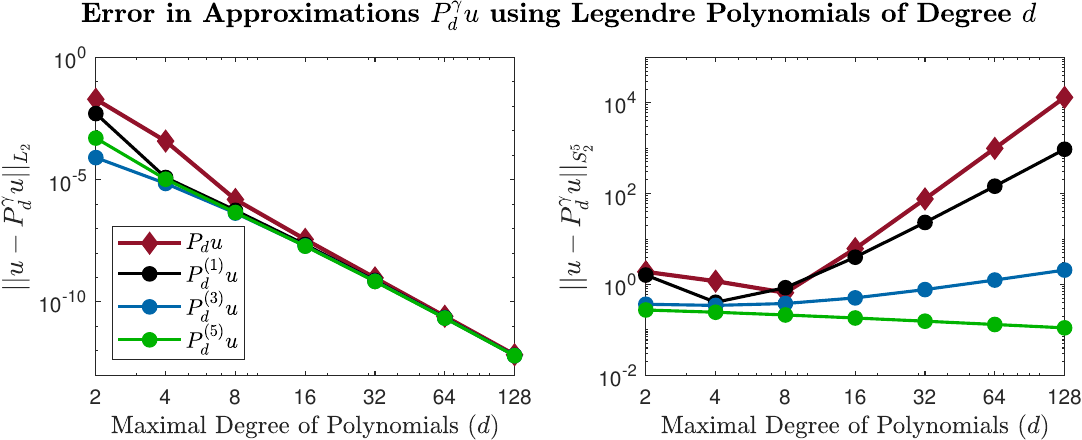}
	\caption{
		$L_{2}$ (left) and Sobolev (right) norms of error in polynomial approximations $P_{d}^{\gamma}u$ of function $u\in S_{2}^{5}[-1,1]$ defined in~\eqref{eq:example1}. 
		The plot $P_{d}u$ corresponds to a standard Legendre projection of $u$ using the $L_{2}$ inner product, as defined in~\eqref{eq:projection_leg}.
		The plots $P_{d}^{\gamma}u$ correspond to reconstructing an approximation of $u$ from a projection of $\partial_{s}^{\gamma}u$ as in~\eqref{eq:projection_leg_sob}, using Theorem~\ref{thm:fundamental_expansion_ND}. 
		Projecting the highest-order derivative $\partial_{s}^{5}u$ of $u\in S_{2}^{5}$, both the $L_{2}$ and Sobolev norms of the error in the associated approximation $P_{d}^{(5)}u$ decrease with the degree $d$, with the $L_{2}$ norm decaying at the same rate as observed for the standard projection $P_{d}u=P_{d}^{(0)}u$.
	}
	\label{fig:invsqrt_conv_Legendre}
\end{figure}

\begin{figure}[t]
	\centering
	\includegraphics[width=1.0\textwidth]{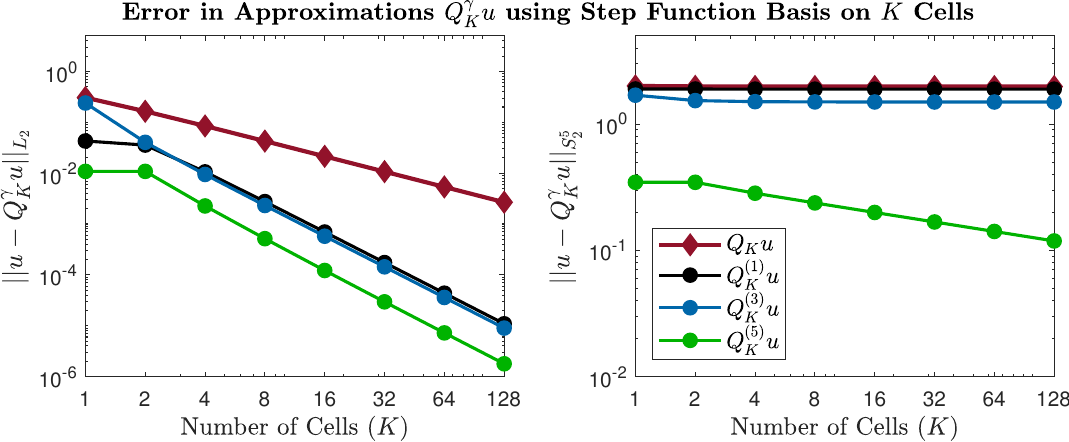}
	\caption{$L_{2}$ (left) and Sobolev (right) norms of error in step function approximations $Q_{K}^{\gamma}u$ of $u\in S_{2}^{5}[-1,1]$ defined in~\eqref{eq:example1}, on a uniform grid of $K$ cells. Each $Q_{K}^{\gamma}u$ is computed by projecting $D^{\gamma}u$ onto a space spanned by $K$ orthogonal step functions, and subsequently reconstructing an approximation of $u$ as in~\eqref{eq:projection_stepfun_sob}. The error in the direct step function projection $Q_{K}u=Q_{K}^{0}u$ is also plotted.}
	\label{fig:invsqrt_conv_stepfun}	
\end{figure}

To illustrate the proposed polynomial and step function approximation methods, each method is applied to construct an approximation of the univariate function $u\in S_{2}^{5}[-1,1]$ defined by
\begin{equation}\label{eq:example1}
	u(s)=\frac{s^4}{36} +\frac{17s^3}{210} -\frac{3s^2}{55} +\frac{29s}{90}-\frac{413}{1140} +\text{sign}(s)\frac{512|s|^{\frac{19}{4}}}{65835}.
\end{equation}
Note that, since $\Omega=[-1,1]$ is only one-dimensional in this case, $S_{2}^{5}[-1,1]=W_{2}^{5}[-1,1]$ by definition of the Sobolev spaces. Since $\partial_{s}^{5}u(s)=\frac{1}{2|s|^{1/4}}$ is square integrable on $[-1,1]$, but $\partial_{s}^{6}u(s)=\frac{-s}{8|s|^{9/4}}$ is not, $u$ is properly an element of $S_{2}^{5}[-1,1]=W_{2}^{5}[-1,1]$.
Figure~\ref{fig:invsqrt_conv_Legendre} and Figure~\ref{fig:invsqrt_conv_stepfun} show the error in the polynomial projection $P_{d}^{\gamma}u$ and step function approximation $Q_{K}^{\gamma}u$ for $\gamma=0,1,3,5$, and for increasing polynomial degrees $d\in\N$ and numbers of cells $K\in\N$, respectively. The error in each approximation is computed both in the $L_{2}$ norm and the standard norm on $S_{2}^{5}$, the latter of which in this case coincides with the standard Sobolev norm on $W_{2}^{5}$.

Figure~\ref{fig:invsqrt_conv_Legendre} shows that the $L_{2}$ norm of the error in the polynomial approximation $P_{d}^{\gamma}u$ decreases at a linear rate on log-log scale, independent of the order $\gamma$ of the Sobolev projection. The slope of each graph is roughly $-5$, matching the expected rate of decay $\|u-P_{d}u\|_{L_{2}}=\mcl{O}(d^{-\delta})=\mcl{O}(d^{-5})$ for the standard $L_{2}$ projection. On the other hand, the error in the Sobolev norm decreases only for $\gamma=5$, displaying a decay $\|u-P_{d}^{(5)}u\|_{S_{2}^{5}}=\mcl{O}(d^{-0.25})$.

Figure~\ref{fig:invsqrt_conv_stepfun} shows that, just as for the polynomial approximation, the $L_{2}$ norm of the error in the step function approximation decreases linearly on log-log scale. Although the observed rate of decay is substantially smaller than using a Legendre polynomial basis, decaying only as $\mcl{O}(K^{-1})$ for the direct approximation $Q_{K}=Q_{K}^{0}$, this rate of decay is improved by projecting a derivative of the function ($\gamma\in\{1,3,5\}$), achieving decay of $\mcl{O}(K^{-2})$. However, since the approximations $Q_{K}^{\gamma}u\in S_{2}^{\gamma}$ are only differentiable up to order $\gamma$, convergence in the Sobolev norm is observed exclusively when approximating the highest-order derivative ($\gamma=5$) of the function, though still achieving decay of only $\mcl{O}(K^{-0.005})$.

\begin{figure}[t]
	\centering
	\includegraphics[width=1.0\textwidth]{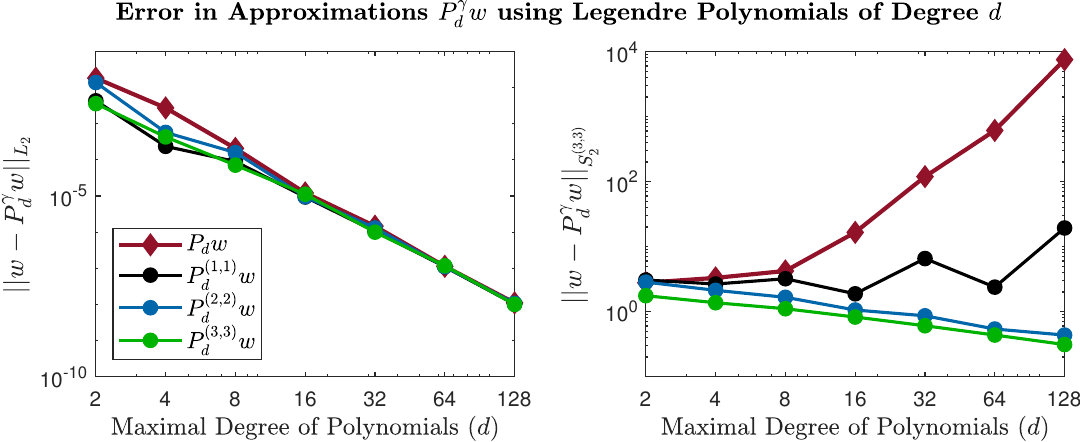}
	\caption{
		$L_{2}$ (left) and Sobolev (right) norms of error in polynomial approximations $P_{d}^{\gamma}w$ of $w\in S_{2}^{(3,3)}[[-1,1]^2]$ defined in~\eqref{eq:example2}.
		Each approximation $P_{d}^{\gamma}w$ is computed by projecting the derivative $D^{\gamma}w$ onto the space of polynomials of degree at most $d$ in each variable, and reconstructing an approximation of $w$ as in~\eqref{eq:projection_leg_sob}. The error in the direct projection $P_{d}w=P_{d}^{(0,0)}w$ is also plotted. 
		Convergence of the approximation in the Sobolev norm is observed only for $P_{d}^{(2,2)}w$ and $P_{d}^{(3,3)}w$, i.e. when projecting $D^{(2,2)}w$ and $D^{(3,3)}w$.
	}
	\label{fig:blockwave_2D_conv_Legendre}
\end{figure}

\subsubsection{Numerical Example in 2D}

\begin{figure}[t]
	\centering
	\includegraphics[width=1.0\textwidth]{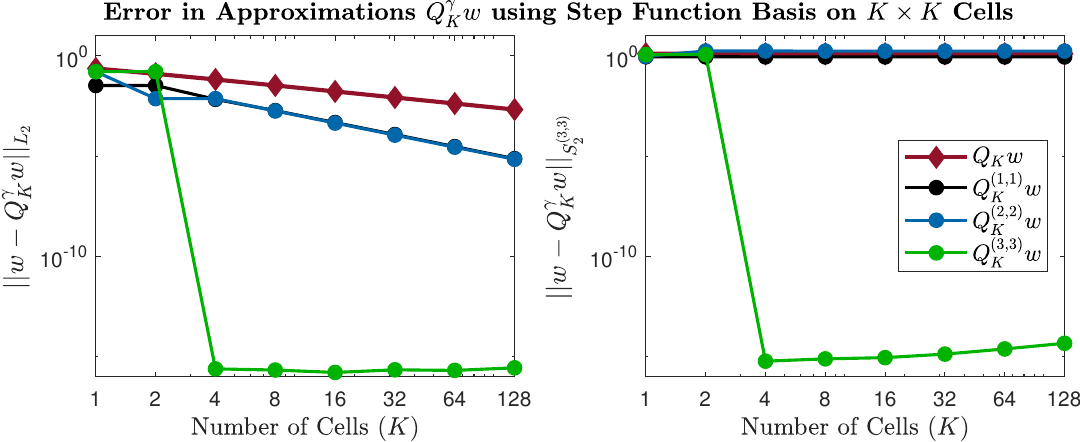}
	\caption{$L_{2}$ (left) and Sobolev (right) norms of error in step function approximations $Q_{K}^{\gamma}w$ of $w\in S_{2}^{(3,3)}[[-1,1]^2]$ defined in~\eqref{eq:example2}, on a grid of $K\times K$ cells. Each $Q_{K}^{\gamma}w$ is computed by projecting $D^{\gamma}w$ onto a space spanned by $K\times K$ orthogonal step functions, and subsequently reconstructing an approximation of $w$ as in~\eqref{eq:projection_stepfun_sob}. The error in the direct step function projection $Q_{K}w=Q_{K}^{(0,0)}w$ is also plotted. Note that $Q_{K}^{(3,3)}w=w$ for $K\geq 4$, since $D^{(3,3)}w$ is a piecewise constant function on $4\times 4$ cells.}
	\label{fig:blockwave_conv_stepfun}	
\end{figure}

For a second example, consider the function $w\in S_{2}^{(3,3)}[[-1,1]^2]$ defined by $w(s_{1},s_{2})=v(s_{1})v(s_{2})$, where for $x\in[-1,1]$,
\begin{equation}\label{eq:example2}
	v(x)
	=\begin{cases}
		\frac{1}{6}x^3 - \frac{1}{2}x^2 + \frac{5}{6}x - \frac{1}{6},	&	x<-\frac{1}{2},	\\
		-\frac{1}{6}x^3 +\frac{7}{12}x - \frac{5}{24},	&	|x|\leq \frac{1}{2},	\\
		\frac{1}{6}x^3 - \frac{1}{2}x^2 + \frac{5}{6}x - \frac{1}{4},	&	x>\frac{1}{2}.
	\end{cases}
\end{equation}
Note that the third-order derivative of $v$ is a step function,
\begin{equation*}
	\partial_{x}^{3}v=
	\begin{cases}
		1,	&	|x|>\frac{1}{2},	\\
		-1,	&	\text{else},
	\end{cases}
\end{equation*}%
which is square-integrable but not differentiable on $[-1,1]$. As such, the weak derivative $D^{\alpha}w(s_{1},s_{2})=\partial_{s_{1}}^{\alpha_{1}}v(s_{1})\partial_{s_{2}}^{\alpha_{2}}v(s_{2})$ is only well-defined for $(0,0)\leq\alpha\leq(3,3)$, and $w$ is properly an element of $S_{2}^{(3,3)}[[-1,1]^2]$.
Figure~\ref{fig:blockwave_2D_conv_Legendre} and Figure~\ref{fig:blockwave_conv_stepfun} show the error in the polynomial projection $P_{d}^{\gamma}w$ and step function approximation $Q_{K}^{\gamma}w$ for $\gamma=(\gamma_{0},\gamma_{0})$ with $\gamma_{0}\in\{0,1,2,3\}$, and for increasing polynomial degrees $(d,d)$ and numbers of cells $K\times K$, respectively. As in the 1D example, the error in each approximation is computed both in the $L_{2}$ norm and the standard norm on $S_{2}^{(3,3)}$, though in this case we remark that $\|\mbf{u}\|_{S_{2}^{(3,3)}}\geq\|\mbf{u}\|_{W_{2}^{3}}$ for any $\mbf{u}\in S_{2}^{(3,3)}$.

The figures show that, again, the $L_{2}$ norm of the error decreases linearly on log-log scale using both projection methods. Using the Legendre polynomial basis, the observed decay is roughly $\mcl{O}(d^{-\delta_{0}})=\mcl{O}(d^{-3})$, independent of the order $\gamma$ of the Sobolev projection $P_{d}^{\gamma}w$. Using the step function basis, this decay is only $\mcl{O}(K^{-1})$ for the direct projection $Q_{K}=Q_{K}^{0}$, and $\mcl{O}(K^{-2})$ when projecting a derivative of order $\gamma_{0}\in\{1,2\}$. However, since the highest-order derivative $D^{(3,3)}w$ of the function $w$ is itself a step function on $4\times 4$ cells, the approximation $Q_{4}^{(3,3)}w$ now exactly matches the function $w$ itself, and both the $L_{2}$ and $S_{2}^{(3,3)}$ norm of the error in this approximation drop to machine precision for $K=4$. For every other value $\gamma_{0}\in\{0,1,2\}$, the approximations $Q_{K}^{\gamma}w$ cannot capture derivatives $D^{\alpha}w$ for $\alpha>\gamma$, and the Sobolev norm of the error does not decay. Using Legendre polynomials, decay of the error in the Sobolev norm is observed only for $\gamma=(2,2)$ and $\gamma=(3,3)$, i.e. when projecting a suitably high order derivative of the function.


\section{Conclusion}

In this paper, a representation of differentiable functions on a hyperrectangle was proposed, expressing such functions in terms of their highest-order derivatives, as well as their lower-order derivatives along lower boundaries of the domain. This representation was proven to offer a bijective map between the dominating mixed smoothness Sobolev space $S_{2}^{\delta}[\Omega]$ and a space of square-integrable functions $\mbf{L}_{2}^{\delta}[\Omega]$, defined by an integral operator with polynomial kernel. Based on this map, a method for polynomial approximation of differentiable functions was proposed, by projecting derivatives of the function onto a space of polynomials or step functions, and then using the main theorem to reconstruct an approximation of the original function. Through numerical examples, it was shown that this approach exhibits similar convergence to standard projection methods in the $L_{2}$ norm, and performs substantially better in the Sobolev norm. While the map proposed in this paper is defined only for functions on a hyperrectangle, similar maps can likely be derived for functions on more general domains, using integral operators with more general (non-polynomial) kernels.

\bibliographystyle{elsarticle-num}
\bibliography{bibfile}


\appendix

\section{Boundary Values of Functions in $S_{2}^{\delta}[\Omega]$}\label{appx:SobolevSpace}


Theorem~\ref{thm:fundamental_expansion_ND} proves that every differentiable function $\mbf{u}\in S_{2}^{\delta}[\Omega]$ on a hyperrectangle $\Omega:=\prod_{i=1}^{N}[a_{i},b_{i}]$ can be expressed in terms of a suitable set of derivatives and boundary values $B^{\alpha-\delta}D^{\alpha}\mbf{u}$, for $0^{N}\leq\alpha\leq\delta$. In this appendix, these boundary values are proven to indeed be well-defined.

To start, we first remark that for arbitrary $\delta\geq 1^{N}$, it has already been proven in~\cite{abdulla2023SobolevEmbedding} that $S_{2}^{\delta}[\Omega]$ can be embedded into a space of H\"older continuous functions, and in fact, that $D^{\alpha}\mbf{u}$ is continuous for every $\mbf{u}\in S_{2}^{\delta}[\Omega]$ and $0^{N}\leq \alpha\leq\delta-1^{N}$. As such, it is clear that the values of these derivatives $D^{\alpha}\mbf{u}$ along the boundaries of the hyperrectangle are well-defined, and the focus in this appendix will instead be on the case that $\delta\leq 1^{N}$. In particular, fix arbitrary $k\in\{1,\hdots,N\}$, and consider $\delta=\tnf{e}_{k}\in\R^{N}$ given by the $k$th standard Euclidean basis vector, so that
\begin{equation*}
	S_{2}^{\delta}[\Omega]=S_{2}^{\tnf{e}_{k}}[\Omega]:=\{\mbf{u}\in L_2[\Omega]\mid \partial_{s_{k}}\mbf{u}\in L_2[\Omega]\}.
\end{equation*}
We will show that any function $\mbf{u}\in S_{2}^{\tnf{e}_{k}}[\Omega]$ is absolutely continuous along $s_{k}\in[a_{k},b_{k}]$, and thus well-defined at $s_{k}=a_{k}$ and $s_{k}=b_{k}$. To this end, let
\begin{equation}\label{eq:Omega_tilde_k}
	\tilde{\Omega}_{k}:=[a_{1},b_{1}]\times\cdots\times[a_{k-1},b_{k-1}]\times[a_{k+1},b_{k+1}]\times\cdots\times[a_{N},b_{N}] \subseteq\R^{N-1}.
\end{equation}
We prove that any function $\mbf{u}\in S_{2}^{\tnf{e}_{k}}[\Omega]$ can be identified with a function $\mbf{v}:[a_{k},b_{k}]\to L_2[\tilde{\Omega}_{k}]$ for which both $\mbf{v}$ and $\partial_{s_{k}}\mbf{v}$ are bounded. Specifically, let $L_2([a_{k},b_{k}],L_2[\tilde{\Omega}_{k}])$ denote the Lebesge-Bochner space of equivalence classes of square-integrable functions mapping $[a_{k},b_{k}]\to L_{2}[\tilde{\Omega}_{k}]$, with the norm
\begin{equation*}
	\|\mbf{v}\|_{L_2}^2=\int_{a_{k}}^{b_{k}}\|\mbf{v}(s_{k})\|_{L_{2}}^2 ds_{k},\qquad \mbf{v}\in L_2([a_{k},b_{k}],L_{2}[\tilde{\Omega}_{k}]).
\end{equation*}
Then, there exists an isometric isomorphism $L_{2}[\Omega]\cong L_{2}([a_{k},b_{k}],L_{2}[\tilde{\Omega}_{k}])$, identifying $\mbf{v}\in L_{2}([a_{k},b_{k}],L_{2}[\tilde{\Omega}_{k}])$ with associated $\mbf{u}\in L_{2}[\Omega]$ through the relation $\mbf{u}(s)=\mbf{v}(s_{k})(\tilde{s})$ for all $s\in \Omega$ and $\tilde{s}\in\tilde{\Omega}_{k}$~\cite{defant1992TensorNorms}. Now, a Sobolev space of functions $\mbf{v}:[a_{k},b_{k}]\to L_{2}[\tilde{\Omega}_{k}]$ can be similarly defined as
\begin{equation*}
	W_{2}^{1}([a_{k},b_{k}],L_{2}[\tilde{\Omega}_{k}])
	:=\{\mbf{u}\in L_2([a,b],L_{2}[\tilde{\Omega}_{k}])\mid \partial_{s_{k}}\mbf{u}\in L_2([a_{k},b_{k}],L_{2}[\tilde{\Omega}_{k}])\}
\end{equation*}
with the norm
\begin{equation*}
	\|\mbf{v}\|_{W_{2}^{1}}^2=
	\int_{a_{k}}^{b_{k}}\|\mbf{v}(s_{k})\|_{L_{2}}^2 ds_{k} +\int_{a_{k}}^{b_{k}}\|\partial_{s_{k}}\mbf{v}(s_{k})\|_{L_{2}}^2 ds_{k},\quad \mbf{v}\in W_2^{1}([a_{k},b_{k}],L_{2}[\tilde{\Omega}_{k}]).
\end{equation*}
Using the isomorphic relation between $L_{2}[\Omega]$ and $L_{2}([a_{k},b_{k}],L_{2}[\tilde{\Omega}_{k}])$, it is not difficult to see that elements in $S_{2}^{\tnf{e}_{k}}[\Omega]$ can be identified with associated elements in $W_{2}^{1}([a_{k},b_{k}],L_{2}[\tilde{\Omega}_{k}])$. In particular, we have the following result.

\begin{prop}\label{prop:SobolevIsomorphism}
	Let $N\in\N$, $\delta\in\N_{0}^{N}$ and $\Omega=\prod_{i=1}^{N}[a_{i},b_{i}]$. Let further $k\in\{1,\hdots,N\}$, and define $\tilde{\Omega}_{k}$ as in~\eqref{eq:Omega_tilde_k}. Then, $S_{2}^{\tnf{e}_{k}}[\Omega] \cong W_{2}^{1}([a_{k},b_{k}],L_{2}[\tilde{\Omega}_{k}])$ isometrically. Specifically, elements $\mbf{u}\in S_{2}^{\tnf{e}_{k}}[\Omega]$ can be identified with elements $\mbf{v}\in W_{2}^{1}([a_{k},b_{k}],L_{2}[\tilde{\Omega}_{k}])$ through the relation
	\begin{equation*}
		\mbf{v}(s_{k})(s_{1},\hdots,s_{k-1},s_{k+1},\hdots,s_{N})=\mbf{u}(s_{1},\hdots,s_{N}),\qquad
		\forall (s_{1},\hdots,s_{N})\in\Omega.
	\end{equation*}
\end{prop}
\begin{pf}
	Without loss of generality, let $k=N$, changing the order of the variables if necessary. Write $\tilde{\Omega}:=\tilde{\Omega}_{N}=\prod_{i=1}^{N-1}[a_{i},b_{i}]$.
	
	To prove that $S_{2}^{\tnf{e}_{N}}[\Omega] \hookrightarrow  W_{2}^{1}([a_{N},b_{N}],L_{2}[\tilde{\Omega}])$, fix arbitrary $\mbf{u}\in S_{2}^{\tnf{e}_{N}}[\Omega]$, and let associated $\mbf{v}$ be such that $\mbf{v}(s_{N})(\tilde{s})=\mbf{u}(\tilde{s},s_{N})$ for all $(\tilde{s},s_{N})\in\tilde{\Omega}\times[a_{N},b_{N}]$. Then $\mbf{u},\partial_{s_{N}}\mbf{u}\in L_2[\Omega]$, and therefore $\mbf{v},\partial_{s_{N}}\mbf{v}\in L_2([a_{N},b_{N}],L_2[\tilde{\Omega}])$. It follows that $\mbf{v}\in W_2^{1}([a_{N},b_{N}],L_{2}[\tilde{\Omega}])$.
	
	Next, to prove that $W_{2}^{1}([a_{N},b_{N}],L_{2}[\tilde{\Omega}]) \hookrightarrow S_{2}^{\tnf{e}_{N}}[\Omega]$, fix arbitrary $\mbf{v}\in W_{2}^{1}([a_{N},b_{N}],L_{2}[\tilde{\Omega}])$, and define associated $\mbf{u}$ by $\mbf{u}(\tilde{s},s_{N})=\mbf{v}(s_{N})(\tilde{s})$ for all $(\tilde{s},s_{N})\in\tilde{\Omega}\times[a_{N},b_{N}]$. Then $\mbf{v},\partial_{s_{N}}\mbf{v}\in L_2([a_{N},b_{N}],L_{2}[\tilde{\Omega}])$, and therefore $\mbf{u},\partial_{s_{N}}\mbf{u}\in L_2[\Omega]$. It follows that $\mbf{u}\in S_{2}^{\tnf{e}_{N}}[\Omega]$.
	
	We find $S_{2}^{\tnf{e}_{N}}[\Omega] \hookrightarrow  W_{2}^{1}([a_{N},b_{N}],L_{2}[\tilde{\Omega}])$ and $W_{2}^{1}([a_{N},b_{N}],L_{2}[\tilde{\Omega}]) \hookrightarrow S_{2}^{\tnf{e}_{N}}[\Omega]$, and therefore $S_{2}^{\tnf{e}_{N}}[\Omega] \cong W_{2}^{1}([a_{N},b_{N}],L_{2}[\tilde{\Omega}])$. To see that this isomorphism is isometric, fix arbitrary $\mbf{v}\in W_{2}^{1}([a_{N},b_{N}],L_{2}[\tilde{\Omega}])$ and associated $\mbf{u}\in S_{2}^{\tnf{e}_{N}}[\Omega]$. Then, by Fubini's theorem,
	\begin{align*}
		\|\mbf{v}\|_{W_{2}^{1}}^2
		&=\int_{a_{N}}^{b_{N}}\|\mbf{v}(s_{N})\|_{L_2}^2 ds_{N} +\int_{a_{N}}^{b_{N}}\|\partial_{s_{N}}\mbf{v}(s_{N})\|_{L_2}^2 ds_{N} \\
		&=\int_{a_{N}}^{b_{N}}\int_{\tilde{\Omega}}|\mbf{v}(s_{N})(\tilde{s})|^2 d\tilde{s}\: ds_{N} +\int_{a_{N}}^{b_{N}}\int_{\tilde{\Omega}}|\partial_{s_{N}}\mbf{v}(s_{N})(\tilde{s})|^2 d\tilde{s}\: ds_{N}	\\
		&=\int_{\Omega}|\mbf{u}(s)|^2 ds +\int_{\Omega}|\partial_{s_{N}}\mbf{u}(s)|^2 ds
		=\|\mbf{u}\|_{S_{2}^{\tnf{e}_{N}}}^{2}.
	\end{align*}
	Hence, $\|\mbf{v}\|_{W_{2}^{1}}=\|\mbf{u}\|_{S_{2}^{\tnf{e}_{N}}}$, and thus the isomorphism is isometric.
\end{pf}

Proposition~\ref{prop:SobolevIsomorphism} proves that any function $\mbf{u}\in S_{2}^{\tnf{e}_{k}}[\Omega]$ is isomorphic to an associated element $\mbf{v}\in W_{2}^{1}([a_{k},b_{k}],L_{2}[\tilde{\Omega}_{k}])$. Using the fact that the Sobolev space $W^{1}_{2}$ on an interval $[a_{k},b_{k}]$ can be embedded into a space of continuous functions on $[a_{k},b_{k}]$, it follows that the boundary elements $\mbf{u}|_{s_{k}=a_{k}}$ and $\mbf{u}|_{s_{k}=b_{k}}$ are well-defined. In particular, we have the following corollary.%

\begin{cor}\label{cor:FTC_Sobolev}
	For $N\in\N$, $k\in\{1,\hdots,N\}$ and $\Omega:=\prod_{i=1}^{N}[a_{i},b_{i}]\subseteq\R^{N}$, if $\mbf{u}\in S_{2}^{\tnf{e}_{k}}[\Omega]$, then 
	\begin{equation*}
		\mbf{u}=\mbf{u}|_{s_{k}=a_{k}} +\int_{a_{k}}^{s_{k}}(\partial_{s_{k}}\mbf{u})|_{s_{k}=\theta}d\theta,\qquad
		\tnf{a.e. }s_{k}\in(a_{k},b_{k}).
	\end{equation*}
\end{cor}
\begin{pf}
	Fix arbitrary $\mbf{u}\in S_{2}^{\tnf{e}_{k}}[\Omega]$, and let $\tilde{\Omega}_{k}$ be as in~\eqref{eq:Omega_tilde_k}. Then, by Proposition~\ref{prop:SobolevIsomorphism}, and by Corollary~1.4.36 in~\cite{cazenave1998SemilinearEqs},
	\begin{align*}
		\mbf{u}\in S_{2}^{\tnf{e}_{k}}[\Omega]
		\cong W_{2}^{1}([a_{k},b_{k}],L_{2}[\tilde{\Omega}_{k}])	
		\hookrightarrow C_{\tnf{b,u}}((a_{k},b_{k}),L_{2}[\tilde{\Omega}_{k}]),
	\end{align*} 
	where $C_{\tnf{b,u}}$ denotes the space of uniformly continuous, bounded functions from $(a_{k},b_{k})$ to $L_{2}[\tilde{\Omega}_{k}]$. It follows that the boundary element $\mbf{u}|_{s_{k}=a_{k}}\in L_{2}[\tilde{\Omega}_{k}]$ is well-defined. Moreover, by Theorem~1.4.35 in~\cite{cazenave1998SemilinearEqs},
	\begin{equation*}
		\mbf{u}|_{s_{k}=t_{1}} -\mbf{u}|_{s_{k}=t_{0}}
		=\int_{t_{0}}^{t_{1}}(\partial_{s_{k}}\mbf{u})|_{s_{k}=\theta}d\theta,\qquad
		\tnf{a.e. }t_{0},t_{1}\in(a_{k},b_{k}).
	\end{equation*}
	Taking the limit as $t_{0}\to a_{k}$, the proposed result follows.
\end{pf}

By this corollary, a boundary operator $B^{d\cdot\tnf{e}_{k}}:S^{\tnf{e}_{k}}[\Omega]\to L_{2}[\tilde{\Omega}_{k}]$ for $d\in\Z$ can be defined such that $B^{d\cdot\tnf{e}_{k}}\mbf{u}=\mbf{u}|_{a_{k}}$ for $d<0$ and $B^{d\cdot\tnf{e}_{k}}\mbf{u}=\mbf{u}|_{b_{k}}$ for $d>0$. Since $S_{2}^{\delta}[\Omega]\subseteq S_{2}^{\tnf{e}_{k}}[\Omega]$ for any $\delta\geq\tnf{e}_{k}$, it follows that the more general boundary operator $B^{\beta}=\prod_{i=1}^{N}B^{\beta_{i}\tnf{e}_{i}}:S_{2}^{\delta}[\Omega]\to L_{2}[\Omega^{\beta}]$, evaluating $B^{\beta}\mbf{u}:=\mbf{u}|_{\Omega^{\beta}}$, is also well-defined for $\mbf{u}\in S_{2}^{\delta}[\Omega]$, so long as $\beta_{i}=0$ whenever $\delta_{i}=0$.
Since, for any $\mbf{u}\in S_{2}^{\delta}[\Omega]$, we have $D^{\alpha}\mbf{u}\in S_{2}^{\delta-\alpha}[\Omega]$ for all $0^{N}\leq\alpha\leq\delta$, it follows that the boundary values $B^{\alpha-\delta}D^{\alpha}\mbf{u}$ in the expansion of $\mbf{u}$ in Theorem~\ref{thm:fundamental_expansion_ND} are indeed well-defined elements of $L_2[\Omega^{\alpha-\delta}]$.

\end{document}